\documentclass[12pt]{siamltex}
\usepackage{color,graphicx,subfig}
\usepackage{amsbsy,amsmath,amssymb}
\usepackage{colonequals} 
\usepackage{bbm} 
\usepackage{enumerate,comment,caption}
\usepackage[top=2.75cm, bottom=3.00cm, left=2.50cm, right=2.50cm]{geometry}
\usepackage[section]{placeins} 
\usepackage{titlesec}
\usepackage{diagbox}
\titlelabel{\thetitle.\quad}

\numberwithin{equation}{section}

\newcommand{\be}{backward-Euler }

\newcommand{\CN}{Crank-Nicolson }
\newcommand{\fsts}{fractional-step $\theta$-scheme }

\newcommand{\indom}[1]{\quad \quad \mbox{in}\; #1}

\newcommand{\bs}[1]{\boldsymbol{#1}}

\title{Fully implicit time-stepping schemes and non-linear solvers for systems of reaction-diffusion equations}
\date{\today}
\author{Anotida Madzvamuse\footnotemark[2] \footnotemark[4] \and Andy H.W. Chung\footnotemark[3]
\footnotemark[4]}

\begin{document}
\maketitle
\renewcommand{\thefootnote}{\fnsymbol{footnote}}
\footnotetext[2]{corresponding author: a.madzvamuse@sussex.ac.uk}
\footnotetext[4]{University of  Sussex, School of Mathematical and Physical Sciences, Department of
Mathematics, University of Sussex, Brighton, BN1 9QH, United Kingdom}

\begin{abstract}
In this article we present robust, efficient and accurate fully implicit time-stepping schemes and
nonlinear solvers for systems of reaction-diffusion equations. The applications of reaction-diffusion
systems is abundant in the literature, from modelling pattern formation in developmental biology to
cancer research, wound healing, tissue and bone regeneration and cell motility. Therefore, it is
crucial that modellers, analysts and biologists are able to solve accurately and efficiently systems of
highly nonlinear parabolic partial differential equations on complex stationary and sometimes continuously evolving domains and surfaces. The main contribution of our paper is the study of fully implicit schemes by use of the Newton method and the Picard iteration applied to the backward Euler, the Crank-Nicolson (and its modifications) and the fractional-step $\theta$ methods. Our results conclude that the fractional-step $\theta$ method  coupled with a single Newton iteration at each timestep is as accurate as the fully adaptive Newton method; and both outperform the Picard iteration. In particular, the results strongly support the observation that a single Newton iteration is sufficient to yield as accurate results as those obtained by use of an adaptive Newton method. This is particularly advantageous when solving highly complex nonlinear partial differential equations on evolving domains and surfaces. To validate our theoretical results, various appropriate numerical
experiments are exhibited on stationary planary domains and in the bulk of stationary surfaces.
\end{abstract}

\section{Introduction}
Since the pioneering work of Turing \cite{Tur1952}, a wide variety of models of reaction-diffusion equations (RDEs) have been proposed as plausible mechanisms for pattern generation processes  \cite{Mur2003}. Turing derived the conditions under which a linearised reaction-diffusion system admits a linearly stable spatially homogeneous steady state in the absence of diffusion and  yet, in the presence of diffusion, it becomes unstable under appropriate conditions to yield a spatially varying inhomogeneous steady state. This process is now well-known as diffusively-driven instability and is of particular interest in developmental biological pattern formation as a means of initiating self organisation from a virtually homogeneous background. Turing patterns were first observed by Castets {\it et. al.} \cite{Castets1990} in a chloride-ionic-malonic-acid reaction. Ouyang and Swinney \cite{Ouyang1991} were the first to observe a Turing instability from a spatially uniform state to a patterned state. Although controversial in a biological context for many years, recent experimental findings suggest this may be a mechanism for the formation of repeated structures in skin organ formation \cite{Sick2006,Maini2006} and zebrafish mesoderm cell fates \cite{Solnica-Krezel2003}. Beyond developmental biology, reaction-diffusion models are widely applied in cell motility, cancer biology, astrophysics, semiconductor physics, ecology, material science, chemistry, financial mathematics and textile engineering \cite{Ascher1995,Atkinson1989,Bangerth2007,Barreira2011,Batchelor2000,BS2003,Castets1990,Dufiet1992,Madzvamuse2003,Mur2003,O2003}.

In many cases, these models comprise of highly nonlinear reaction terms which makes it impossible to obtain analytical solutions in closed form. Hence numerical methods are employed. For the reaction-diffusion theory of pattern formation, these methods play two key roles: (i) they are used to validate linear stability theoretical results close to bifurcation points (and vice versa) and (ii) far away from bifurcation points, they provide numerical approximate solutions (in the absence of analytical solutions).  A typical numerical method consists of two processes: first, a space discretisation is applied to render the system of partial differential equations (PDEs) into a system of ordinary differentials equations (ODEs) and second, a time discretisation is employed, thereby transforming the system of ODEs into a system of linear or nonlinear algebraic
 equations depending on the form of the time discretisation scheme. Finally, techniques from numerical linear algebra are employed to solve efficiently the resulting system.

 Space discretisation methods include (but are not limited to) finite differences, finite elements, spectral methods, finite volume, closest point methods on stationary domains, volumes and surfaces and more recently moving grid, surface finite element and particle methods applied to domains and surfaces that evolve in time \cite{Madzvamuse2003,Barreira2011}. In this article, we choose to use the finite element methodology since it can cope easily with complex irregular geometries, volumes and surfaces. Other numerical methods can be applied, however, for some methods (e.g., the finite differences) their extension to complicated, irregular and sometimes continuously evolving domains and surfaces is not at all trivial. The applicability of the finite element methods to complicated domains is well known. The finite element method can easily deal with complicated and sometimes continuously changing domains. 

 Several time discretisation schemes have been widely used to compute solutions of PDEs on stationary and evolving domains and surfaces \cite{Mad2006}. These include the forward Euler's method (the most commonly used in computational biology), Gear's method, a modified Euler predictor-corrector method, Gourlay's method \cite{Dufiet1992}, a semi-implicit Rosenbock integrator \cite{Dufiet1992}, and Runge-Kutta schemes \cite{Hairer1987}. Most of these are inadequate because of the stiffness of the diffusive term. Fully explicit methods require excessively small time-steps resulting in computations that are prohibitively expensive especially  in multi-dimensions. Recently, IMEX schemes have been used to solve RDEs
on stationary one-dimensional domains \cite{Ascher1995,Ruuth1995}. The key essence of these schemes is that an implicit scheme is applied to approximate the diffusive term and an explicit scheme is used to approximate the reaction kinetics, hence the name IMEX. More recently, Madzvamuse \cite{Mad2006} presented several time-stepping schemes to compute solutions to reaction-diffusion systems on fixed and growing domains. A first order semi-implicit backward Euler differentiation formula (1-SBEM) which treats diffusive and linear reaction terms implicitly and nonlinear reaction terms semi-implicitly was introduced and shown to be more robust than IMEX schemes. The 1-SBEM employed a single Picard iteration. In all these studies, very little work has been done to extend such analysis to fully implicit schemes. Therefore the focus of this article is to  substantially extend previous IMEX schemes to fully implicit time-stepping schemes for RDEs on stationary domains and surfaces. Fully implicit time-stepping schemes offer greater numerical stability than the 1-SBEM and IMEX schemes. One revealing result of our studies is that, for RDEs, a single Newton iteration is sufficient to obtain as accurate solutions as when an adaptive Newton scheme is used. This single Newton iteration for the nonlinear reaction kinetics, outperforms a single Picard iteration applied to IMEX time-stepping schemes used in \cite{Mad2006}. On the other hand, the overall elapsed CPU time taken by the Picard iteration is substantially lower than that taken by the Newton method for the backward Euler and Crank-Nicholson and its modifications. The only exception is the fractional $\theta$-method where the Newton method outperforms the Picard iteration in terms of the overall elapsed CPU time. This is attributed to the fact that the Newton method uses the GMRES solver to invert the resulting non-symmetric matrices, whereas the Picard iteration uses the CG solver to invert block diagonal symmetric matrices. Our recommendation is that when solving RDEs on stationary and evolving domains and surfaces using fully implicit schemes, it is sufficient to employ single Newton iteration with the fractional $\theta$-method in order to obtain as accurate solutions  as those obtained when an adaptive Newton method is used. This is particularly relevant for problems posed on evolving domains and surfaces \cite{Elliott2012,Lakkis2013}. Numerical experiments demonstrate that the fractional $\theta$-method coupled with a single Newton method is about 130 times faster than the backward Euler and 15 times faster than the Crank-Nicholson and its modifications.  

 Hence, our article is structured as follows: in Section \ref{sec:prelim} we present the model equations under study and the choice of parameter values to be used. The space and time discretisation methods and schemes are detailed in Section \ref{sec:num_methods}. To validate our theoretical predictions, we present experimental order of convergence results in Section \ref{sec:EOC}. Various numerical experiments are presented and discussed in Section \ref{sec:num_results}. In Section \ref{sec:comparison_with_imex_schemes} we demonstrate how fully implicit schemes and appropriate nonlinear solvers outperform standard IMEX schemes. To demonstrate the applicability and generality of the fully implicit scheme, in Section \ref{sec:geometries} we present various solutions of RDEs on stationary planary domains and in the bulk of stationary surfaces.  Finally, we conclude and outline future research in Section \ref{sec:conclusion}

\section{The model equations}
\label{sec:prelim}
 \par Let $\Omega$ be a convex domain with Lipschitz boundary, $I\!=\![0,T]$ be some time interval and
$(\bs{x},t)\in\Omega\times I$. We take for illustrative purposes the well-known Schnakenberg reaction kinetics (also known as the {\it activator-depleted} substrate or Brusselator model \cite{GM1972,PL1968,Sch1979})
 to obtain the model system of RDEs which reads
 \begin{equation} 
 \left\{
 \begin{aligned}
 & \frac{ \partial u}{\partial t} - \nabla^2u  = \gamma(a-u+u^2v) := \gamma f(u,v), \\
 & \frac{ \partial v}{\partial t} - d\nabla^2v = \gamma(b-u^2v) := \gamma g(u,v),\indom{\Omega\times I},
 \end{aligned}
 \right.  \label{schnak}
 \end{equation}
 for the pair $(u(\bs{x},t),v(\bs{x},t))$ and some real positive numbers $a$, $b$, $d$ and $\gamma$. Here, $u$
and $v$ correspond to concentrations of  two chemical species, $d$ measures the ratio of the relative diffusivity of the $v$ to $u$ and
$\gamma$ measures the strength of the reaction. To close the system, we choose homogeneous Neumann boundary
conditions on the entire boundary and take initial conditions to be small random
perturbations around the steady state values
 \begin{equation}
 (u_{eq}\,,v_{eq})=\left(a+b\,,\,\frac{b}{(a+b)^2}\right).  \label{eqpoints}
 \end{equation}
Next, we recollect briefly conditions for diffusion-driven instability \cite{Tur1952} used to guide parameter estimation as well as yielding linear approximate solutions close to bifurcation points which validate numerical results.  With no spatial variation in $u$ or $v$, in the absence of diffusion, the equilibrium point
(\ref{eqpoints}) is linearly stable provided that \cite{Mur2003}
 \begin{equation} 
 f_u +  g_v < 0 \quad
 \mbox{and} \quad f_ug_v-f_vg_u > 0 , \label{linstab}
 \end{equation}
where the derivatives are evaluated at the equilibrium point. If one then allows spatial inhomogeneity,
it is possible that the system evolves to an inhomogeneous steady state. This entails an initial
instability caused by diffusion, a phenomenon known as diffusion-driven instability, or Turing
instability after the author who first described it in \cite{Tur1952}. Subsequently, the non-linear
reaction terms keep the solution bounded in an invariant set \cite{Smo1994}.
\par A linear stability analysis reveals conditions that will drive a Turing instability \cite{Mur2003}.
Let us consider small perturbations from the equilibrium and write them as
$\tilde{u}(\bs{x},t)\!=\!u(\bs{x},t)-u_{eq}$ and $\tilde{v}(\bs{x},t)\!=\!v(\bs{x},t)-v_{eq}$. Define
$\bs{\xi}(\bs{x},t)\!=\!(\tilde{u}(\bs{x},t),\tilde{v}(\bs{x},t))$. Then linearising (\ref{schnak}) we
obtain the system of linear PDEs
\begin{equation}
\frac{\partial \bs{\xi} }{\partial t} = \gamma \left(
     \begin{array}{cc}
        f_u & f_v \\
        g_u & g_v \\
         \end{array} \right) \bs{\xi} + \left(
     \begin{array}{cc}
        1 & 0 \\
        0 & d\\
         \end{array} \right) \nabla^2\bs{\xi},
\end{equation}
with homogeneous Neumann boundary conditions.
 This system can be solved by separation of variables to yield
\begin{equation} 
\bs{\xi}(\bs{x},t)=\sum_k \bs{c_k}e^{\lambda t}\xi_k, \label{solunmodedecomp}
\end{equation}
where $\xi_k$ are the modes which solve the homogeneous Neumann problem
\begin{equation} 
\nabla^2{\xi_k}+k^2\xi_k=0. \label{laplaceeigen}
\end{equation}
 These modes will decay with time unless their wavenumber $k$ lies in the range
\begin{equation} 
k_{-}^2< k^2 < k_{+}^{2} , \label{k_range}
\end{equation}
where
\begin{equation} 
k_{\pm}=\gamma\,\frac{df_u+g_v \pm \sqrt{\left( df_u+g_v\right)^2 - 4d(f_ug_v-f_vg_u)}}{2d} ,
\end{equation}
and the following conditions hold:
\begin{equation}
df_u+g_v>0 \quad \mbox{and} \quad (df_u+g_v)^2 - 4d(f_ug_v-f_vg_u)>0,
\end{equation}
where the derivatives are evaluated at the equilibrium values. Thus, if we perturb the system from
equilibrium, under certain choices of parameters, we can expect exponential growth of some modes which
correspond to the linearly unstable modes of (\ref{solunmodedecomp}). This restriction of parameters
defines the so-called Turing space \cite{Tur1952,Mur2003}.

\par When $\Omega$ is the two-dimensional unit square, the eigenmodes of (\ref{laplaceeigen}) have the
form $\cos(n\pi x)\cos(m\pi y)$ for  $n,m\in\mathbb{Z}$ with eigenvalues $k^2\!=\!n^2+m^2$.  The choices
$a\!=\!0.1$, $b\!=\!0.9$, $d\!=\!10$ and $\gamma\!=\!29$ will lead to diffusion-driven instability
\cite{Mad2000}. The parameters chosen guarantee that the modes corresponding to $n^2+m^2=1 $ are
linearly unstable. Furthermore, it can be calculated that the corresponding exponential growth rate is
about $1.6246$ \cite{Mur2003}. This exponential growth rate will later serve as an aid to see if our
numerical simulations are consistent with the theory.

\section{Numerical methods}
\label{sec:num_methods}
\par Due to the non-linearites, an analytical solution to (\ref{schnak}) is not readily available, and
so we seek numerical solutions using the finite element method as detailed below.

\subsection{Finite element spatial discretisation}
To discretise in space we use the finite element method (FEM)
\cite{Hug2003,Tho2006}. We shall experiment with a few time discretisation methods. Time-stepping
schemes for this reaction system were investigated in \cite{Mad2006} for both stationary and evolving
domains. The treatment of the non-linear terms was semi-implicit - here we shall treat them implicitly,
thereby granting greater stability.

\subsection{Space discretisation}
  \par Let $T_h$ be a triangulation of $\Omega_h\subseteq\Omega$, and let each node have coordinates
$\bs{x}_i$, $i=1,2,\cdots,N_h$, with $N_h$ denoting the number of degrees of freedom. Let $\{\phi_i\}$ be the set of piecewise linear shape functions on $T_h$. Now, implementation of the FEM yields \cite{Hug2003,Tho2006}
  \begin{equation} 
 \left\{
 \begin{aligned}
 & M\dot{\bs{u}}+\gamma M \bs{u} + A\bs{u} - \gamma
B(\bs{u},\bs{v})\bs{u}=\gamma a\bs{\mathbbm{1}}_{\phi},  \\
  & \quad M\dot{\bs{v}} + d A\bs{v} + \gamma
B(\bs{u},\bs{u})\bs{v}=\gamma b\bs{\mathbbm{1}}_{\phi},
 \end{aligned}
 \right.       \label{schnakfem}
 \end{equation}
 where $A$ and $M$ are the stiffness and mass matrices respectively with entries
 \begin{equation}
 a_{ij} =\int_{\Omega}{\nabla\phi_{i}\!\cdot\!\nabla\phi_{j}\,d\bs{x}} \quad \mbox{and} \quad
m_{ij}=\int_{\Omega}{\phi_i\phi_j\,d\bs{x}},
 \end{equation}
  $\bs{\mathbbm{1}}_{\phi}$ is the column vector with $j$-th entry $\phi_j$. Given some vectors
$\bs{a}$ and $\bs{b}$, $B(\bs{a},\bs{b})$ is the matrix with entries
\begin{equation} 
B_{ij}  = \sum\limits_{k=1}^{N_h}\sum\limits_{l=1}^{N_h}
a_kb_l\int_{\Omega}{\phi_k\phi_l\phi_i\phi_j\,d\bs{x}} .
\end{equation}
It is readily checked that given a
third vector $\bs{c}$, the matrix $B$ satisfies
\begin{equation} 
B(\bs{a},\bs{b})\bs{c}=B(\bs{a},\bs{c})\bs{b}=B(\bs{c},\bs{b})\bs{a}. \label{cyclicb}
\end{equation}
 \par Equation (\ref{schnakfem}) does not yet lend itself to a numerical solution. First, it is still
continuous with respect to time. Second, the non-linear matrix $B$ does not allow a solution to be
gained in a straightforward manner.

\subsection{Time discretisation}
For illustrative purposes, we will consider the \be (BE) and \CN (CN) methods and the \fsts (FSTS) \cite{Glo2003} with a uniform timestep $\tau$. To proceed, define
\begin{align}
  \bs{G}_1(\bs{u},\bs{v})  &= A\bs{u} +\gamma M\bs{u} - \gamma B(\bs{u},\bs{v})\bs{u}-\gamma
a\bs{\mathbbm{1}}_{\phi} \\
\mbox{and} \quad \bs{G}_2(\bs{u},\bs{v})  &= dA\bs{v} + \gamma B(\bs{u},\bs{u})\bs{v}-\gamma
b\bs{\mathbbm{1}}_{\phi}\,.
\end{align}
 For the BE discretisation we solve at the ($n\!+\!1$)-th timestep
 \begin{equation} 
 \left\{
 \begin{aligned}
 &  M \frac{\bs{u}^{n+1}-\bs{u}^{n}}{\tau} + \bs{G}_1(\bs{u}^{n+1},\bs{v}^{n+1})=0, \\
  &M  \frac{\bs{v}^{n+1}-\bs{v}^{n}}{\tau} + \bs{G}_2(\bs{u}^{n+1},\bs{v}^{n+1})=0.
 \end{aligned}
 \right.       \label{schnakfembe}
 \end{equation}
For the CN we solve
  \begin{equation} 
 \left\{
 \begin{aligned}
 &  M \frac{\bs{u}^{n+1}-\bs{u}^{n}}{\tau} + \frac{1}{2}
\left[\bs{G}_1(\bs{u}^{n+1},\bs{v}^{n+1}+\bs{G}_1(\bs{u}^n,\bs{v}^n)\right]=0, \\
  & M \frac{\bs{v}^{n+1}-\bs{v}^{n}}{\tau} + \frac{1}{2}
\left[\bs{G}_2(\bs{u}^{n+1},\bs{v}^{n+1})+\bs{G}_2(\bs{u}^n,\bs{v}^n)\right]=0.
 \end{aligned}
 \right.     \label{schnakfemcn}
 \end{equation}
  \par For the FSTS, we split the differential operators, i.e. $\bs{G_1}$ and $\bs{G_2}$, into two parts.
In this problem, there is the natural dichotomy of diffusion and reaction; however, here, we shall
simply make the distinction between terms that are linear and those that are not. Once we have this
separation, one splits the timestep into three unequal portions alternately treating each dividend
implicitly and explicitly. For some $\theta\!\in\!(0,1/2)$, we first solve for the pair $(\bs{u}^{n+\theta},\bs{v}^{n+\theta})$
 \begin{equation} 
 \left\{
 \begin{aligned}
 & M \frac{\bs{u}^{n+\theta}-\bs{u}^{n}}{\theta\tau} + A\bs{u}^{n+\theta} +\gamma M\bs{u}^{n+\theta}
=\gamma a\bs{\mathbbm{1}}_{\phi}+ \gamma B(\bs{u}^n,\bs{v}^n)\bs{u}^n, \\
  & \quad M \frac{\bs{v}^{n+\theta}-\bs{v}^{n}}{\theta\tau} + d A\bs{v}^{n+\theta}  =\gamma
b\bs{\mathbbm{1}}_{\phi}- \gamma B(\bs{u}^n,\bs{u}^n)\bs{v}^n,
 \end{aligned}
 \right.   \label{schnakfefstseta1}
 \end{equation}
 then for the pair $(\bs{u}^{n+1-\theta},\bs{v}^{n+1-\theta})$ we solve
 \begin{equation} 
 \left\{
 \begin{aligned}
 & M \frac{\bs{u}^{n+1-\theta}-\bs{u}^{n+\theta}}{(1-2\theta)\tau} -\gamma
B(\bs{u}^{n+1-\theta},\bs{v}^{n+1-\theta})\bs{u}^{n+1-\theta} =\gamma a\bs{\mathbbm{1}}_{\phi} -
A\bs{u}^{n+\theta} - \gamma M\bs{u}^{n+\theta},  \\
  & M \frac{\bs{v}^{n+1-\theta}-\bs{v}^{n+\theta}}{(1-2\theta)\tau} + \gamma
B(\bs{u}^{n+1-\theta},\bs{u}^{n+1-\theta})\bs{v}^{n+1-\theta} =\gamma b\bs{\mathbbm{1}}_{\phi} - d
A\bs{v}^{n+\theta},
 \end{aligned}
 \right.       \label{schnakfefstseta2}
 \end{equation}
 and finally for the pair $(\bs{u}^{n+1},\bs{v}^{n+1})$ we solve
  \begin{equation} 
 \left\{
 \begin{aligned}
 & M  \frac{\bs{u}^{n+1}-\bs{u}^{n+1-\theta}}{\theta\tau} + A\bs{u}^{n+1} +\gamma M\bs{u}^{n+1} =\gamma
a\bs{\mathbbm{1}}_{\phi}+ \gamma B(\bs{u}^{n+1-\theta},\bs{v}^{n+1-\theta})\bs{u}^{n+1-\theta}, \\
  & \quad M \frac{\bs{v}^{n+1}-\bs{v}^{n+1-\theta}}{\theta\tau} + d A\bs{v}^{n+1}  =\gamma
b\bs{\mathbbm{1}}_{\phi}- \gamma B(\bs{u}^{n+1-\theta},\bs{u}^{n+1-\theta})\bs{v}^{n+1-\theta}.
 \end{aligned}
 \right. \label{schnakfefstseta3}
 \end{equation}
In the first and third step, the non-linear terms are treated explicitly and the linear terms implicitly.
Conversely, in the second step, the non-linear terms are treated implicitly and the linear terms
explicitly. We take $\theta\!=\!1-1/\sqrt{2}$ since it can be shown that this method is second
order convergent in time for this value of $\theta$ \cite{Glo2003}.

\subsection{Techniques for treating the non-linearities}

\par To deal with the non-linearities, we shall use the Picard iteration and the Newton method \cite{QSS2006}.
\par At the ($n\!+\!1$)-th timestep, the ($k\!+\!1$)-th Picard iterate of the BE system (\ref{schnakfembe}) is the solution of the matrix system
\begin{flalign}
\renewcommand*{\arraystretch}{1.5}
 \left( \begin{array}{cc}
  \left(\frac{1}{\tau}+\gamma\right)M+A -\gamma B(\bs{u}^{n+1}_k,\bs{v}_k^{n+1}) & \bs{0} \\
   \bs{0} & \frac{1}{\tau}\,M + dA + \gamma B(\bs{u}_k^{n+1},\bs{u}^{n+1}_k) \end{array}\right)
  \left( \begin{array}{c}
\bs{u}^{n+1}_{k+1} \\
\bs{v}^{n+1}_{k+1}
   \end{array}\right) \nonumber \\[.3cm]
 && \renewcommand*{\arraystretch}{1.5}
   \mkern-36mu \mkern-36mu \mkern-36mu \mkern-36mu \mkern-36mu
    = \left( \begin{array}{c}
 \gamma a\bs{\mathbbm{1}}_{\phi} +\frac{1}{\tau}\,M\bs{u}^{n}  \\
 \gamma b\bs{\mathbbm{1}}_{\phi} + \frac{1}{\tau}\,M\bs{v}^{n}
  \end{array} \right).
  \label{picard_iter}
 \end{flalign}
We note here that the IMEX scheme considered in \cite{Mad2006} is equivalent to taking a single Picard iteration per timestep. Now, take the left-hand side (LHS) of (\ref{schnakfembe}) and define
 \begin{align}  
    &\bs{F}_1(\bs{u}^{n+1},\bs{v}^{n+1})  = \left[ \left(\frac{1}{\tau}+\gamma\right)M+A-\gamma
B(\bs{u}^{n+1},\bs{v}^{n+1})\right]\bs{u}^{n+1}  -\gamma a\bs{\mathbbm{1}}_{\phi}
-\frac{1}{\tau}\,M\bs{u}^{n}, \\
     &\bs{F}_2(\bs{u}^{n+1},\bs{v}^{n+1})  = \left[ \frac{1}{\tau}\,M +dA + \gamma
B(\bs{u}^{n+1},\bs{u}^{n+1}) \right] \bs{v}^{n+1} -\gamma b\bs{\mathbbm{1}}_{\phi}
-\frac{1}{\tau}\,M\bs{v}^{n},
 \end{align}
and let $\bs{F}\!=\!(\bs{F}_1,\bs{F}_2)$. Then we wish to solve
$\bs{F}(\bs{u}^{n+1},\bs{v}^{n+1})\!=\!0$.  The ($k\!+\!1$)-th Newton iterate of the BE system is the solution of
\begin{equation}  
  J_{\bs{F}}|_{(\bs{u}^{n+1}_k,\bs{v}^{n+1}_k)}
   \Big(\bs{u}^{n+1}_{k+1}-\bs{u}^{n+1}_k,\bs{v}^{n+1}_{k+1}-\bs{v}^{n+1}_k\Big)  =
-\bs{F}(\bs{u}^{n+1}_k,\bs{v}^{n+1}_k),  \quad i=1,2, \label{schnaknewt0}
 \end{equation}
where $J_{\bs{F}}$ is the Jacobian matrix of $\bs{F}$ and the column vector
$\Big(\bs{u}^{n+1}_{k+1}-\bs{u}^{n+1}_k,\bs{v}^{n+1}_{k+1}-\bs{v}^{n+1}_k\Big)$ is understood to be
the vertical concatenation of the first and second arguements. Using property (\ref{cyclicb}) of the
matrix $B$, we have for some vector $\bs{\xi}\!\in\!\mathbb{R}^{N_h}$ the G\^{a}teaux derivative
\begin{equation}
\frac{\partial\bs{F}_1(\bs{u}^{n+1}_k,\bs{v}^{n+1}_k)}{\partial\bs{u}^{n+1}}\bs{\xi} = \left[ \Big( \frac{1}{\tau}+\gamma\Big)M+A-2\gamma
B(\bs{u}^{n+1}_{k},\bs{v}^{n+1}_{k})\right]\bs{\xi}.
 \end{equation}
 Similarly we have,
 \begin{align}
 \mbox{} &\frac{\partial\bs{F}_1(\bs{u}^{n+1}_k,\bs{v}^{n+1}_k)}{\partial\bs{v}^{n+1}}\bs{\xi}=
 -\gamma B(\bs{u}^{n+1}_{k},\bs{u}_k^{n+1})\bs{\xi}, \\
 \mbox{} &\frac{\partial\bs{F}_2(\bs{u}^{n+1}_k,\bs{v}^{n+1}_k)}{\partial\bs{u}^{n+1}}\bs{\xi}=
 2\gamma B(\bs{u}^{n+1}_{k},\bs{v}_k^{n+1})\bs{\xi}, \\
 \mbox{and} \quad &
\frac{\partial\bs{F}_2(\bs{u}^{n+1}_k,\bs{v}^{n+1}_k)}{\partial\bs{v}^{n+1}}\bs{\xi}=
 \left[ \frac{1}{\tau}M+dA+\gamma
B(\bs{u}^{n+1}_{k},\bs{u}^{n+1}_{k})\right]\bs{\xi}.
  \end{align}
Thus, at the ($n\!+\!1)$-th  timestep, the ($k\!+\!1)$-th Newton iterate (\ref{schnaknewt0}) is the
solution of the system
\begin{flalign} 
\renewcommand*{\arraystretch}{1.5}
 \left( \begin{array}{cc}
  \left(\frac{1}{\tau}+\gamma\right)M + A - 2\gamma
       B(\bs{u}^{n+1}_k,\bs{v}^{n+1}_k) & - \gamma
B(\bs{u}^{n+1}_k,\bs{u}^{n+1}_k) \\
   2\gamma
B(\bs{u}^{n+1}_k,\bs{v}^{n+1}_k) & \frac{1}{\tau}\,M +  dA + \gamma
B(\bs{u}^{n+1}_k,\bs{u}^{n+1}_k) \end{array}\right)
  \left( \begin{array}{c}
\bs{u}^{n+1}_{k+1} \\
\bs{v}^{n+1}_{k+1}
   \end{array}\right) \nonumber  \\[.3cm]
 && \renewcommand*{\arraystretch}{1.5}
  \mkern-36mu \mkern-36mu \mkern-36mu \mkern-36mu \mkern-36mu \mkern-36mu \mkern-36mu \mkern-36mu
   = \left( \begin{array}{c}
 \gamma a\bs{\mathbbm{1}}_{\phi} +\frac{1}{\tau}\,M\bs{u}^{n} -2\gamma B(\bs{u}^{n+1}_k,\bs{v}^{n+1}_k)\bs{u}^{n+1}_{k}  \\
 \gamma b\bs{\mathbbm{1}}_{\phi} + \frac{1}{\tau}\,M\bs{v}^{n} +2\gamma B(\bs{u}^{n+1}_k,\bs{u}^{n+1}_k)\bs{v}^{n+1}_k
  \end{array} \right).
  \label{newton_iter}
  \end{flalign}
We can easily show that the Jacobian is not singular in our case. Indeed, suppose that
$J_{\bs{F}}|_{(\bs{u}^{n+1}_k,\bs{v}^{n+1}_k)}\bs{\eta}\!=\!\bs{0}$ for some vector
$\bs{\eta}\!\in\!\mathbb{R}^{2N_h}$. Writing $\bs{\eta}\!=\!(\bs{\eta}_1,\bs{\eta}_2)$, with
$\bs{\eta}_1$, $\bs{\eta}_2\!\in\!\mathbb{R}^{N_h}$, we can add the two resulting equations and
rearrange to obtain
\begin{equation}
\begin{aligned}
\bs{\eta}_1
&=-\left[\left(\frac{1}{\tau}+\gamma\right)M+A\right]^{-1}\left(\frac{1}{\tau}M+dA\right)\bs{\eta}_2
\\
           &\colonequals -K_1^{-1}K_2\,\bs{\eta}_2, \\
\end{aligned}
\end{equation}
where we have used the fact that $A$ and $M$ are both positive definite matrices and are, therefore,
invertible. Upon substitution, this then implies that
\begin{equation}
K_2^{-1} \Big( 2\gamma B(\bs{u}^{n+1}_k,\bs{v}^{n+1}_k) K_1^{-1} K_2 - \gamma
B(\bs{u}^{n+1}_k,\bs{u}^{n+1}_k)\Big)\bs{\eta}_2=\bs{\eta}_2.
\label{imposseqn}
\end{equation}
If $\bs{\eta}_2\!\neq\!=\!\bs{0}$, then the matrix on the LHS of (\ref{imposseqn}) must be the identity
matrix. However, this cannot be since the coefficients $\bs{u}^{n+1}_k$ and $\bs{v}^{n+1}_k$ appear on
the LHS but, obviously, do not appear in the identity matrix. Thus, $\bs{\eta}_2$ must be the zero
vector. Similarly, $\bs{\eta}_1$ must also be the zero vector, from which we conclude that $J_{\bs{F}}$
is invertible. A similar argument will show that the relevant matrix for Picard iteration is also
invertible, though it does not indicate whether or not the iterations converge to the desired solution.
For this we need to show that we have a contraction mapping (for proof see the appendix).
\par We can expect linear convergence with the Picard iteration and quadratic convergence with the Newton
method \cite{QSS2006}. The Newton method may fail to converge if the initial guess is far from the
solution \cite{QSS2006}. However, in our case if we have a small enough timestep, we expect that the
solution $(u^n,v^n)$ will provide an adequate initial guess to find the solution $(u^{n+1},v^{n+1})$
since the change in the solution through successive timesteps will be small. Evidently, the matrix to
be inverted in the Picard system (\ref{picard_iter}) is symmetric whereas for the Newton system (\ref{newton_iter})
it is not. This then precludes the use of efficient iterative techniques such as the conjugate gradient method
to find a solution to the Newton system. Instead, for the Newton method, we use GMRES in order to solve the non-symmetric system of linear algebraic equations \cite{Saad1992}.

\section{Experimental order of convergence (EOC)}
\label{sec:EOC}
There is no known analytical solution to the Schnakenberg system (\ref{schnak}). Therefore, we cannot check directly how well our numerical experiments compare with the exact solution. In the following, we use a well known method of constructing a solution that will satisfy a modified version of (\ref{schnak}) from which we can indirectly gauge the performance of subsequent discretisation. Though this is no substitute for error estimation, it is, nevertheless, a useful guide. The simulations in this section and the next were carried out using the software {\it deal.II} \cite{Bangerth2007}.
\par In the following, we take the domain $\Omega$ to be the unit square. Define
\begin{equation} 
\Xi(x,y,t)=\left(\frac{x^3}{3}-\frac{x^2}{2}\right)\left(\frac{y^3}{3}-\frac{y^2
}{2}\right)\left(1+e^{-t}\right).
\end{equation}	
Then $u=v=\Xi$ is the exact solution to the modified system
 \begin{equation} 
  \left\{
 \begin{aligned}
 & \frac{ \partial u}{\partial t} - \nabla^2u - \gamma(-u+u^2v) = \Xi_t
-\nabla^2 \Xi -\gamma(-\Xi+\Xi^3), \\
 & \frac{ \partial v}{\partial t} - d\nabla^2v + \gamma u^2v = \Xi_t
-d \nabla^2 \Xi +\gamma\,\Xi^3, \indom{\Omega}\,,\mbox{for } t>0\,,
 \end{aligned}
 \right.  \label{schnakmod}
 \end{equation}
with homogeneous Neumann boundary conditions. Initial conditions are then defined by 
\begin{equation}
u_0=v_0=2\left(\frac{x^3}{3}-\frac{x^2}{2}\right)\left(\frac{y^3}{3}-\frac{y^2}{
2}\right).
\end{equation}
Here we shall also take the parameter values $a\!=\!0.1$, $b\!=\!0.9$, $d\!=\!10$ and $\gamma\!=\!29$ \cite{Mad2000}. It is easily seen that $u$ and $v$ both tend to the inhomogeneous steady $u_0/2$ as $t\!\rightarrow\!\infty$.
\par If we now solve (\ref{schnakmod}) using the proposed time-step methods we can calculate the error from the exact solution $\Xi$ at each timestep. This is done using five different timesteps $\tau_i=2^{-i}$, $i=1,2,3,4,5$ in the time interval $t\in[0,10]$. In the following, the convergence criteria for the Picard iteration and the Newton method was chosen to be $||u^{n+1}_{k+1}-u^{n+1}_k||<10^{-5}$, and similarly for the variable $v$, where $||\!\cdot\!||$ denotes the $L_2$-norm. A measure of the error from the exact solution at the ($n+1$)-th timestep is given by
\begin{equation}
e_{u}(t^{n+1})=||u^{n+1}-\Xi(t^{n+1})||,
\end{equation}
and a measure of the error of the whole simulation is given by
\begin{equation}
E_{u}= e_{u}(10),
\end{equation}
which is the error at the end time. Similar quantities can be defined for the $v$-variable. We have the error estimate \cite{Tho2006}
\begin{equation}
E_u \leq C[u]\left(h^2+\tau^{\alpha}\right), \label{parabfemerror}
\end{equation}
where $C[u]$ is a constant and $h$ is the greatest length of the squares in the mesh $T_h$. For a particular time-stepping scheme, the maximal value of $\alpha$ such that the above holds is then its order of convergence. Thus, for BE we expect $\alpha\!=\!1$ and for the CN and the FSTS we expect $\alpha\!=\!2$.
\par Let $E_{u,i}$ ($i\!=\!1,2,\cdots,5$) denote the errors found as above using the same time-stepping scheme with time step $\tau_i$. The estimate (\ref{parabfemerror}) is a sharp estimate so that if $h^2\sim\tau^{\alpha}$, we have approximately
\begin{equation} 
\alpha_i\approx\frac{\log(E_{u,i})-\log(E_{u,i-1})}{\log(\tau_i)-\log(\tau_{i-1})}, \quad i>1,
\end{equation}
with the $\alpha_i$ converging to the expected value as the errors decrease. For the second order methods, the coupling $h_i\!=\!\tau_i$ was obtained by constructing an $n\times n$ square grid with $n\!=\!2^i$. For the backward-Euler method $h_i\!=\!\sqrt{\tau_i}$ was taken by constructing a square grid as before but taking the nearest $n$ such that the value of  $1/n$ was closest to $\sqrt{\tau_i}$.

\begin{table}[t]
      \centering
      \footnotesize
        \begin{tabular}{ l | c c c }
                    & BE & CN & FSTS \\ \hline 
         $\alpha_1$ & -1.965$\times10^{-5}$ & 2.049 & 2.064 \\ 
         $\alpha_2$ & 1.217 & 2.020 & 2.020 \\
         $\alpha_3$ & 0.847 & 2.006 & 2.006 \\
         $\alpha_4$ & 1.184 & 2.002 & 2.002 \\
         $\alpha_5$ & 0.835 & -     &  - \\
         $\alpha_6$ & 0.922 & -     &  - \\
         $\alpha_7$ & 1.083 & -     &  - \\
         $\alpha_8$ & 1.048 & -     &  - \\
         $\alpha_9$ & 0.953 & -     &  - 
         \end{tabular}
           \caption{Experimental order of convergence (EOC) for the
$u$-variable. Further EOC values were computed for BE to investigate the
convergence to the expected value of $1$. }
           \label{eoctable}
\end{table}
\par
\begin{figure}[]
	\centering
	\subfloat[][]{
\includegraphics[keepaspectratio=true,width=.48\textwidth]{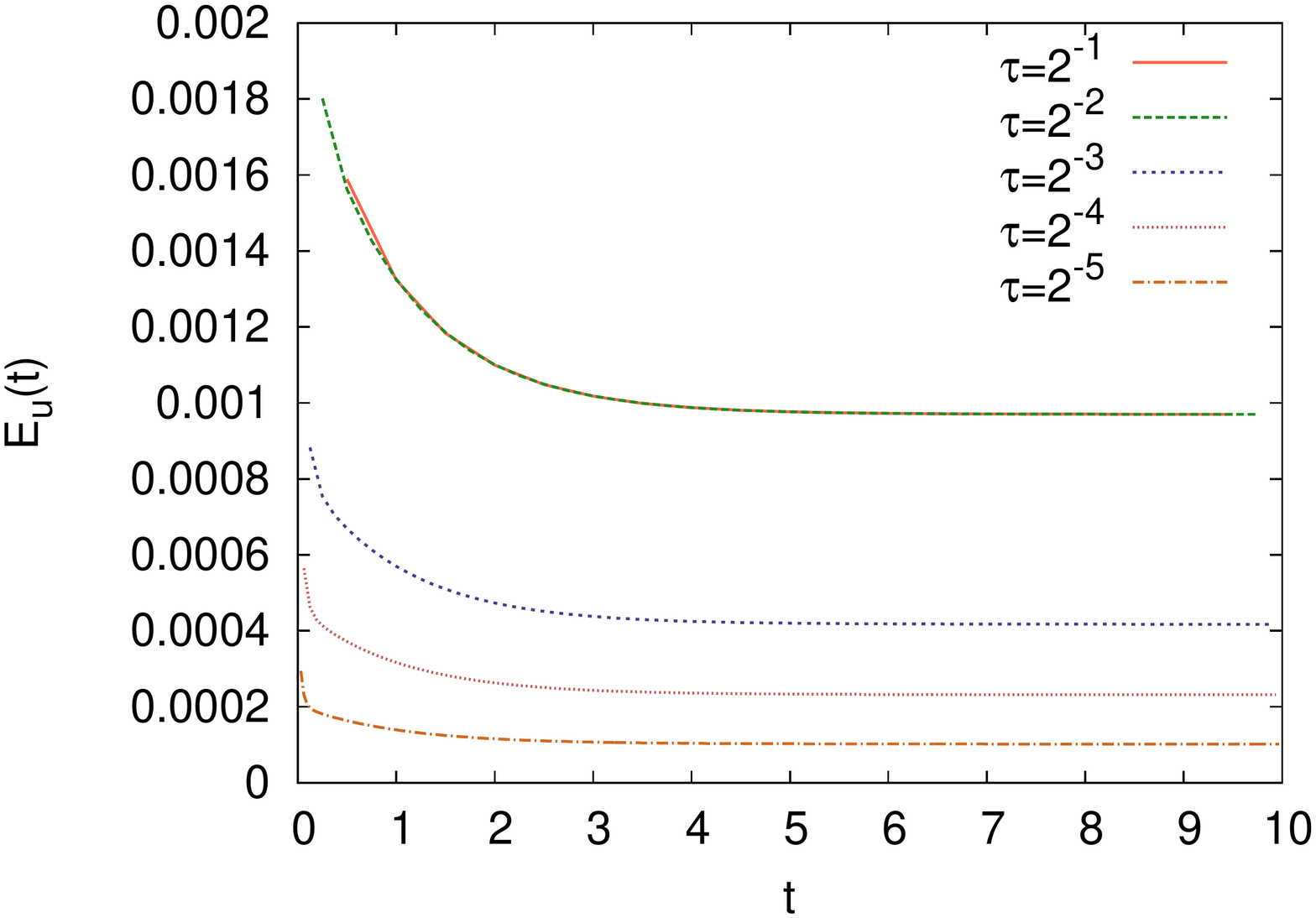} }
	\subfloat[][]{
\includegraphics[keepaspectratio=true,width=0.48\textwidth]{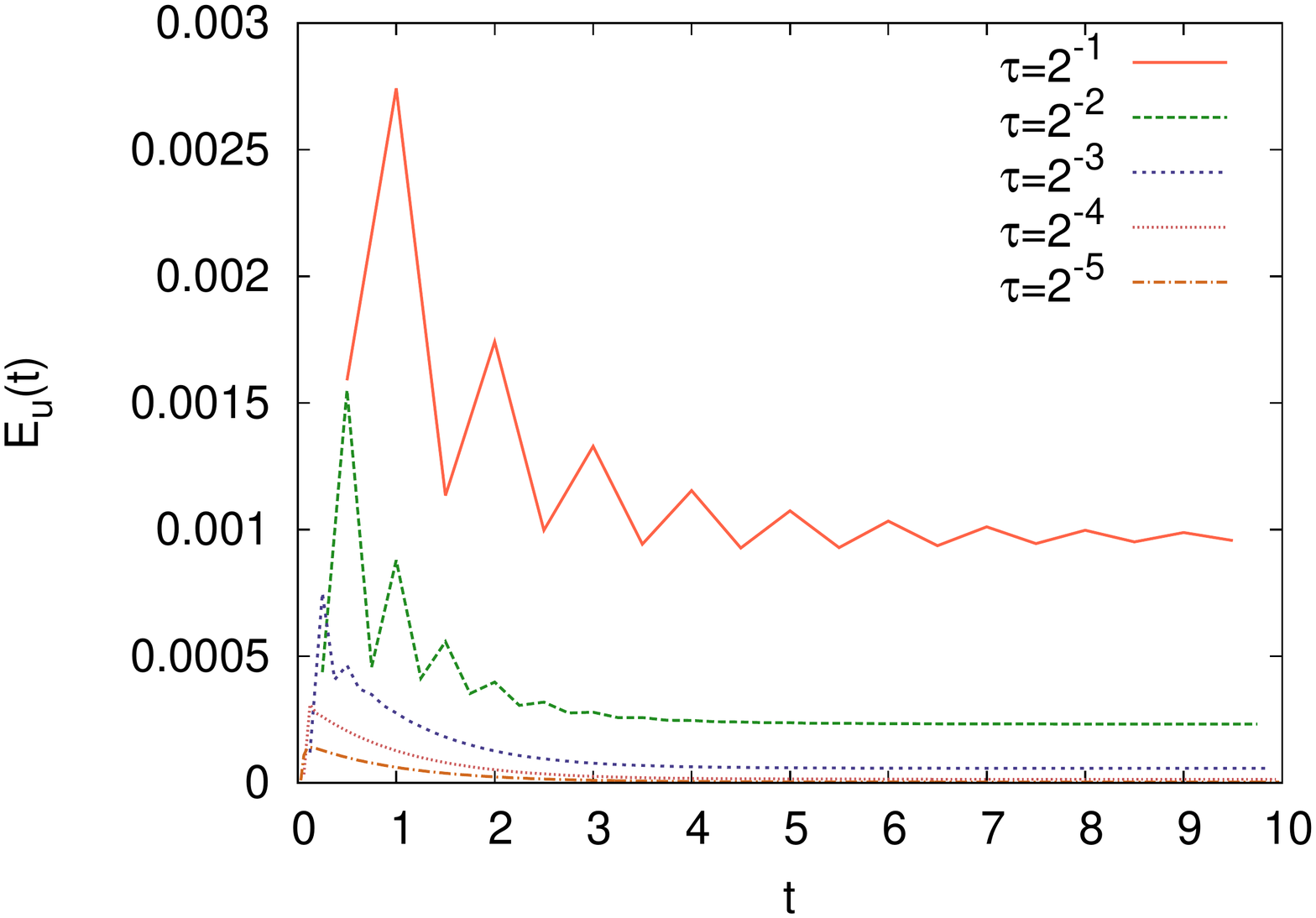} }
     \newline
	\subfloat[][]{
\includegraphics[keepaspectratio=true,width=0.5\textwidth]{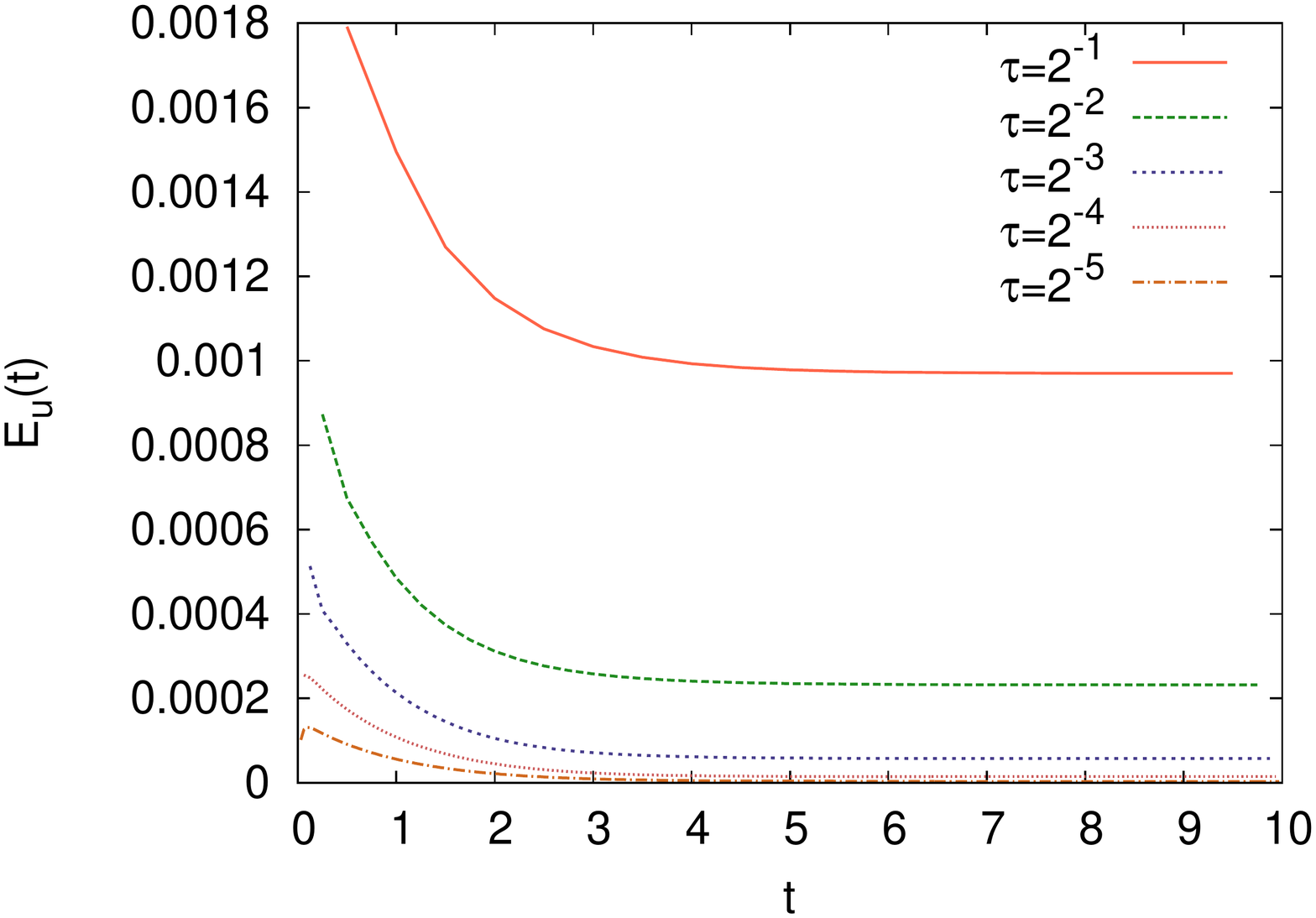} }
\caption{Errors for the modified equation for the $u$-variable using the adaptive Newton method: (a) BE (b) CN (c) FTS.  The second order methods CN and FSTS are more accurate than BE. The error ofconvergence results are shown in Table \ref{eoctable}.}
\label{fig:errors}
\end{figure}

 The errors for $\tau_5=1/32$ are plotted in Fig. \ref{fig:errors}. As expected, the second order methods are more accurate than BE. The convergence results are shown in Table \ref{eoctable}. The values of $\alpha$ obtained are in keeping with the theory.  The apparent non-convergence of the $\alpha_i$ to
$1$ for BE, due perhaps to the errors not being small enough, prompted further values of $\alpha_i$ to be calculated for BE, keeping the timestep $\tau_i\!=\!2^{-i}$. These are shown in Table \ref{eoctable}, where the convergence is made more apparent.

\begin{figure}[]
	\centering
	\subfloat[][]{
\includegraphics[keepaspectratio=true,width=.48\textwidth]{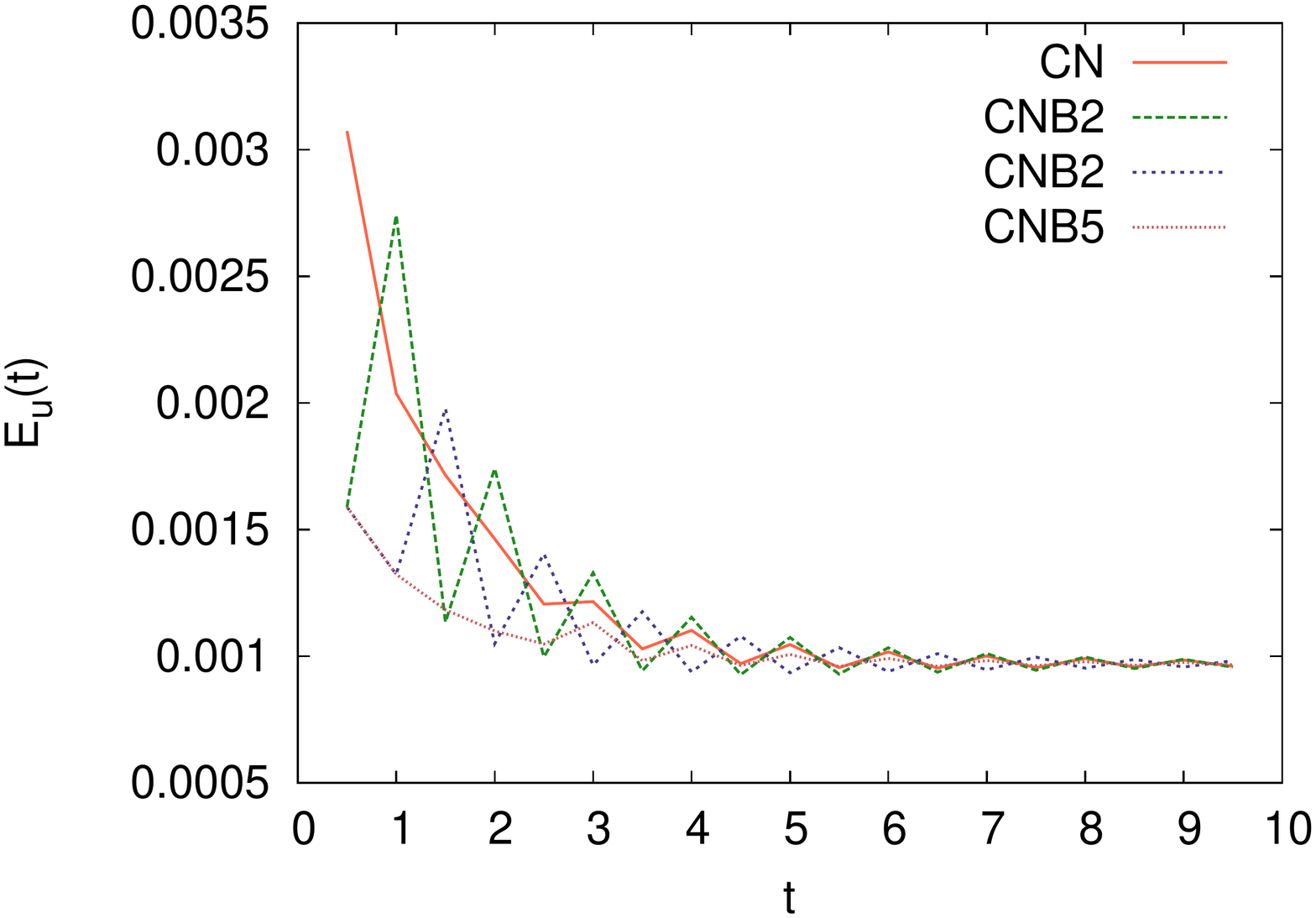}
	\label{fig:ccc1}}
	\subfloat[][]{
\includegraphics[keepaspectratio=true,width=.48\textwidth]{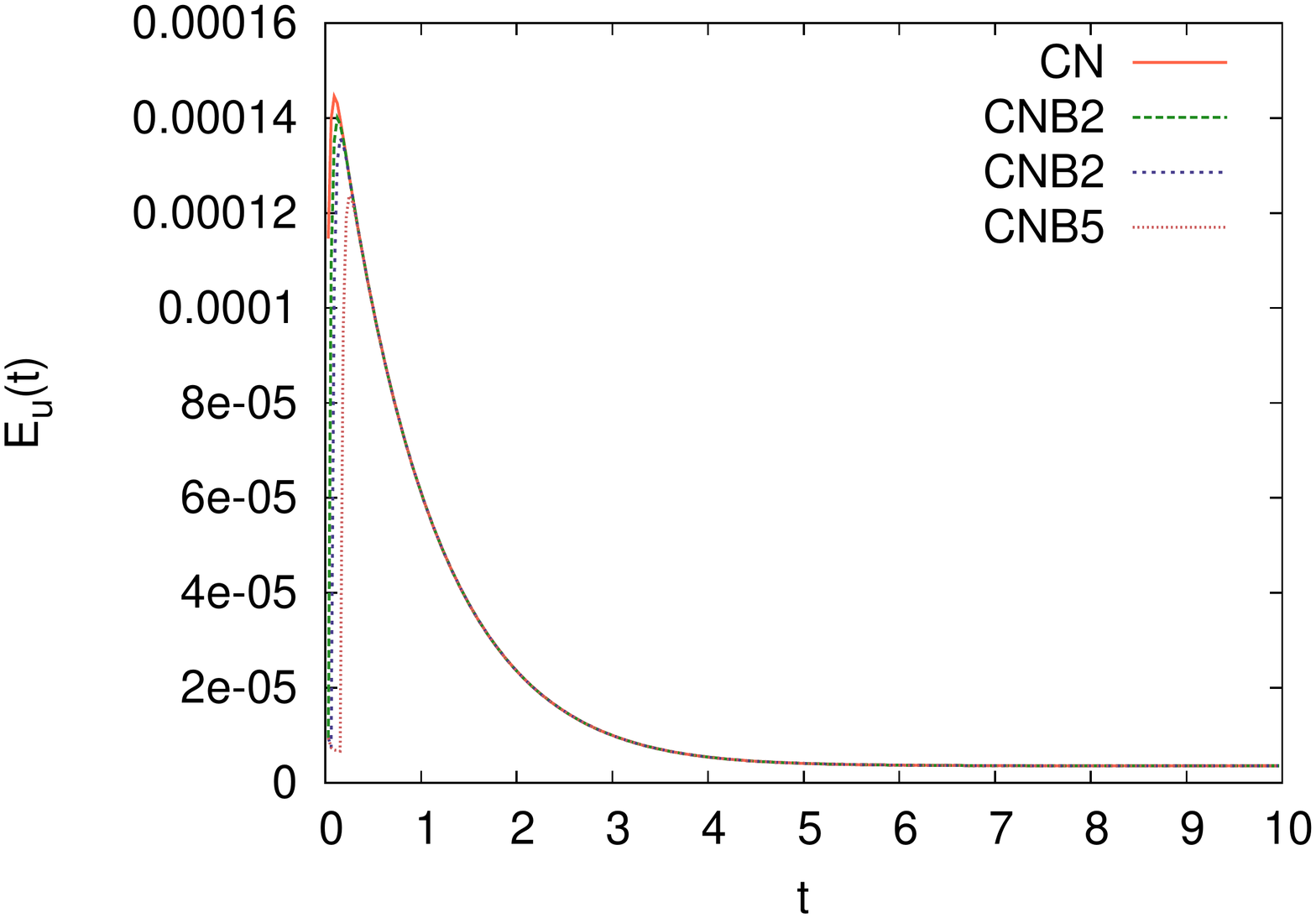}
	\label{fig:ccc2}}
	\caption{Errors in the simulation of the modified equation (\ref{schnakmod}) using the modified
CN methods CNB, CNB2 and CNB5 with (a) $\tau\!=\!1/2$ and (b) $\tau\!=\!1/32$ (see Section \ref{sec:EOC} for details). The oscillations observed in CN are damped by the modified methods for both these timesteps. An adaptive Newton method was employed for the fully implicit time-stepping scheme.}
\label{fig:ccc}
\end{figure}
\par Upon inspection of Fig. \ref{fig:errors}, one can observe oscillations in the error for the CN method which are damped with time. This behaviour is well known \cite{Atkinson1989,Hug2003}, and ultimately mars the results. Put briefly, under the CN, high frequency modes of the errors in the initial data do not damp out effectively \cite{O2003}. Some strategies to reduce these oscillations are studied in \cite{BS2003} and \cite{O2003}. Here, we use a strategy in which BE is performed for a few initial timesteps, afterwhich the CN is implemented. Essentially, the BE method damps out the oscillations which are otherwise retained in the CN method \cite{BS2003}. We try three variants of this method: in the first, we perform BE for the first timestep (CNB); in the second we perform two BE steps (CNB2); and in the third we perform five (CNB5). The results of these implementations are shown in Fig. \ref{fig:ccc} for two timesteps. These methods seem to have varying results. For $\tau\!=\!1/2$, CNB seems to increase the oscillations, whilst CNB5 reduces them the most, however, for the smaller timestep $\tau\!=\!1/32$, all methods seems to behave similarly. Overall however, the magnitude of the errors are reduced taking more BE timesteps initially. 

\section{Numerical Results}
\label{sec:num_results}
We have three time-stepping schemes and two ways to solve the non-linearities implicitly. We choose $\Omega$ to be the unit square and use a uniform square mesh with $10,000$ squares such that the greatest length is $h\!=\!1/\!100$. To ensure errors of a similar magnitude, a timestep of $\tau\!=\!10^{-2}$ was used for the CN method and the FSTS and $\tau\!=\!10^{-4}$ was used for the BE method. For the initial conditions, a random perturbation from equilibrium of order $\sim\!10^{-2}$ was set on every vertex of the mesh. As before, the convergence criterion for the Picard iteration and the Newton method was chosen to be $||u^{n+1}_{k+1}-u^{n+1}_k||<10^{-5}$, and similarly for the variable $v$.
 \par We expect the system to reach a spatially inhomogeneous steady state. Thus, the simulation was left to run until a maximum time $t\!=\!30$, or until the following criterion were satisfied:
 \begin{equation} 
 \frac{||u^{n+1}-u^n||}{\tau}  \leq 10^{-4} \quad \mbox{and} \quad
\frac{||v^{n+1}-v^n||}{\tau} \leq
10^{-4}.  \label{stopcrit}
 \end{equation}
 This quantity is related to the rate of change of the variables, and so we simply stopped the calculations when the solutions have stopped varying significantly with time. We find it convenient to divide by the timestep, $\tau$, to allow for better comparison between solutions using different timesteps.
\begin{table}
      \centering
       \small
      \subfloat[][]{
        \begin{tabular}{ l | c c c c}
                   & BE & CN & CNB5 & FSTS \\ \hline
         end time          & 18.45   & 30.00 & 18.26 & 17.91  \\
         total iterations  &  201024   & 7156  & 4353  & 3822   \\
         elapsed CPU time ($\times 10^4$s) & 185.10 & 25.92 & 10.45 &  11.42\\  
             \label{iterpic}
         \end{tabular}
         }
         \!\!\!\!\!\!\!\!\!\!\!\!\! {\color{white}[.2cm]}
            \subfloat[][]{
        \begin{tabular}{l | c c c c}
                   &  BE & CN & CNB5 & FSTS \\ \hline
         end time         & 18.46   & 30.00  & 18.48 & 18.14 \\
         total iterations & 201052  & 5998   & 3529  & 2994 \\
         elapsed CPU time ($\times 10^4$s) & 963.6  & 105.33 & 96.8  & 8.39
             \label{iternewt}
         \end{tabular}
         }
         \caption{Details of the simulation of the Schnakenberg system using (a)
Picard iteration and (b) the Newton method. Shown are the total number of non-linear iterations required to reach the end time and the elapsed CPU time. Note the non-convergence of the CN method.}
         \label{table:iter}
\end{table}
\par
\begin{figure}[htb]
	\centering
	\subfloat[][]{
\includegraphics[keepaspectratio=true,width=.48\textwidth]{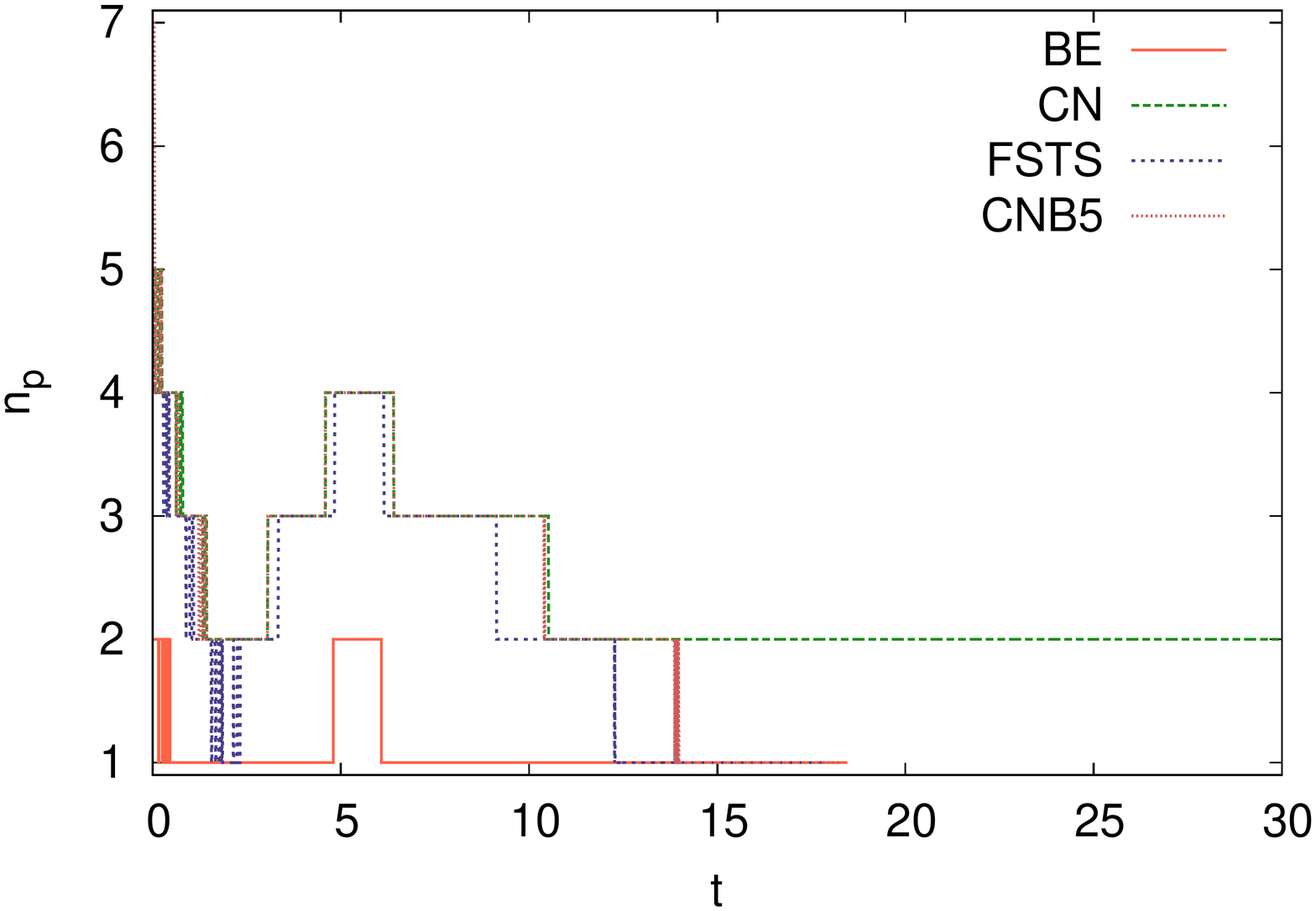}
	\label{iter:p}}
	\subfloat[][]{
\includegraphics[keepaspectratio=true,width=.48\textwidth]{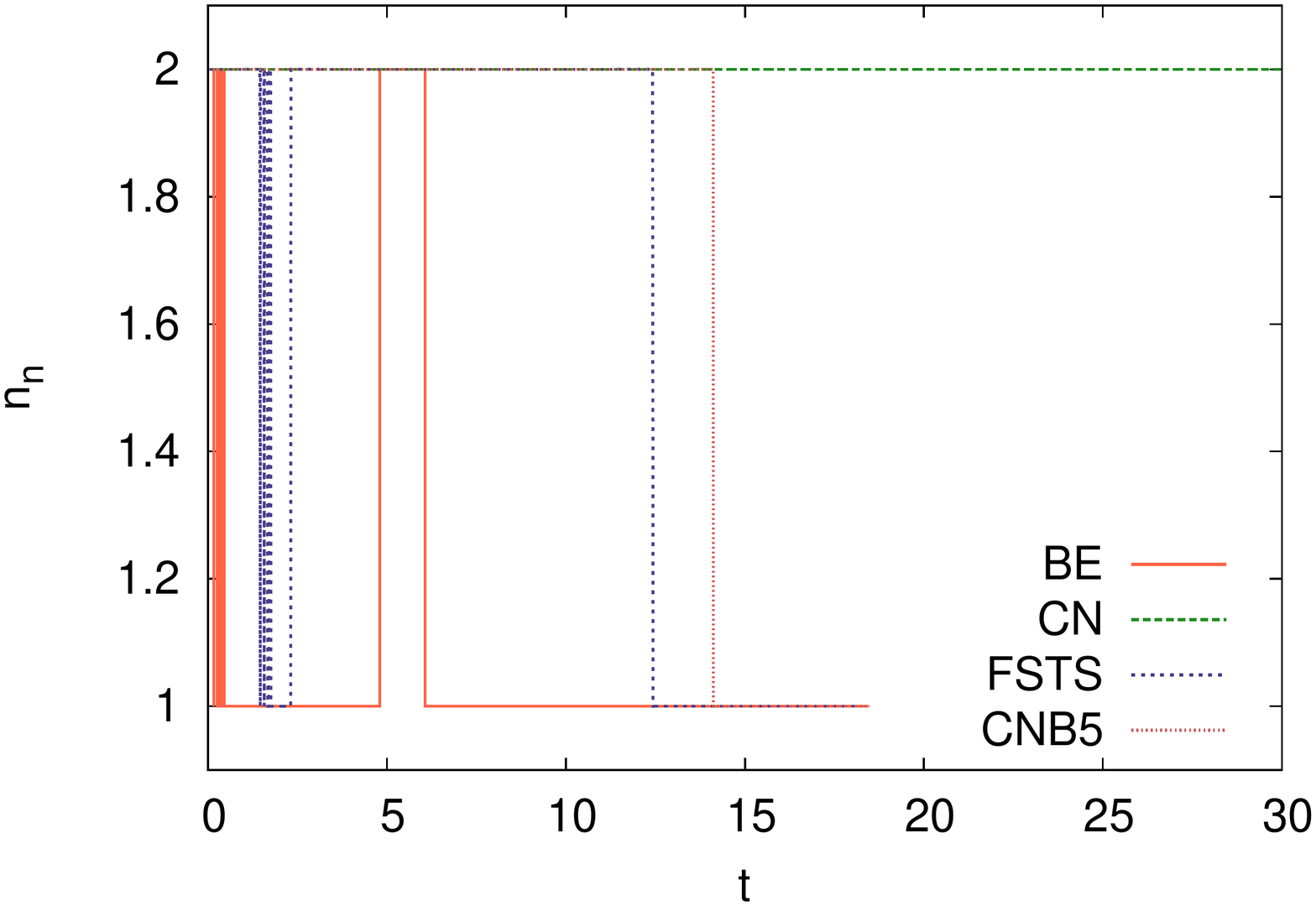}
	\label{iter:n}}
	\caption{Number of non-linear iterations for (a) Picard iteration, $n_p$, and (b) the Newton method, $n_n$, at each timestep.}
	\label{iter}
	\end{figure}

 After performing the calculations, we found that the solutions obtained using the Picard iteration and the Newton method were virtually indistinguishable. Table \ref{table:iter} shows the number of iterations needed to reach the steady state as defined by (\ref{stopcrit}). For the values quoted for FSTS, the two linear steps which frame the non-linear step in each timestep are counted. In Fig. \ref{iter} the number of non-linear iterations performed at each timestep is plotted. The number of Picard iterations varied to a greater degree than those resulting from the Newton method. The number of iterations for the Picard iteration varied between $1$ and $7$, whereas for the Newton method it varied between $1$ and $2$. 

\begin{figure}[tbh]
	\centering
	\subfloat[][]{
\includegraphics[keepaspectratio=true,width=.48\textwidth]{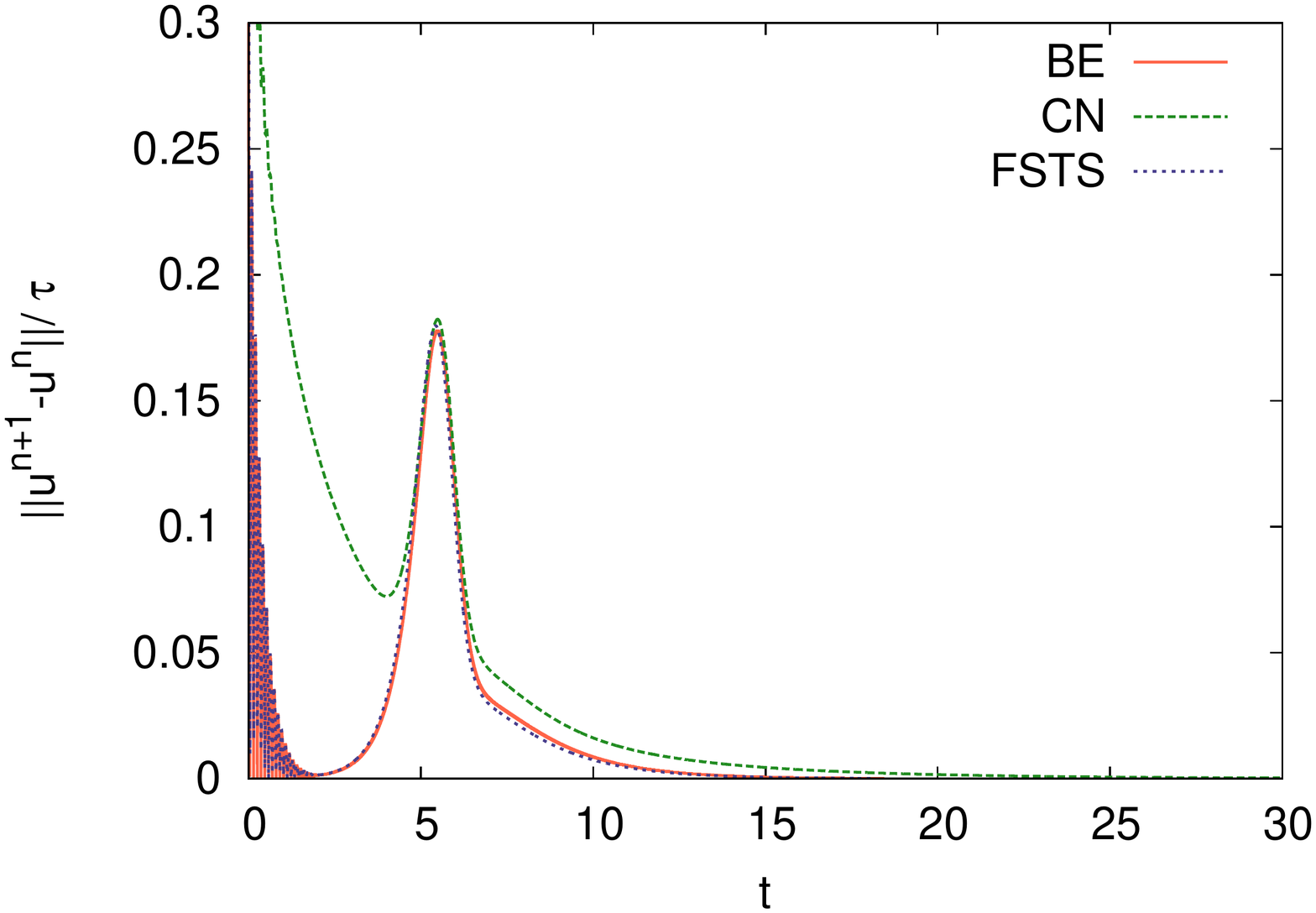}
	\label{bct:u}}
	\subfloat[][]{
\includegraphics[keepaspectratio=true,width=.48\textwidth]{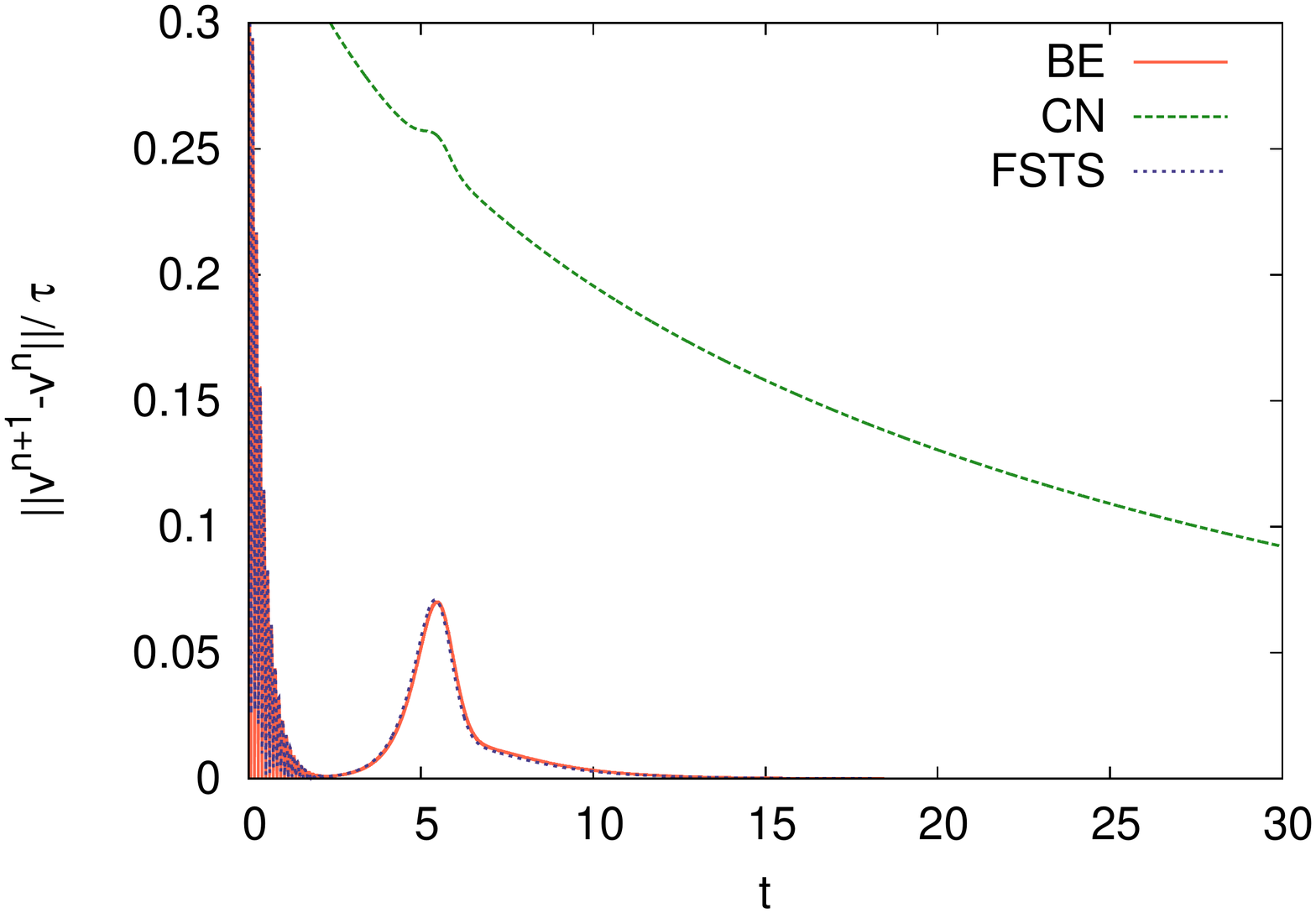}
	\label{bct:v}}
	\caption{Convergence history of the simulation of the Schnakenberg system for (a) the $u$-variable and (b) the $v$-variable. Seen is the initial decay of the modes, mode excitation (growth) and decay into the inhomogeneous steady state. An adaptive Newton method was employed for the fully implicit time-stepping scheme.}
	\label{bct}
	\end{figure}

	

\begin{figure}[tbh]
	\centering
\includegraphics[keepaspectratio=true,width=.6\textwidth]{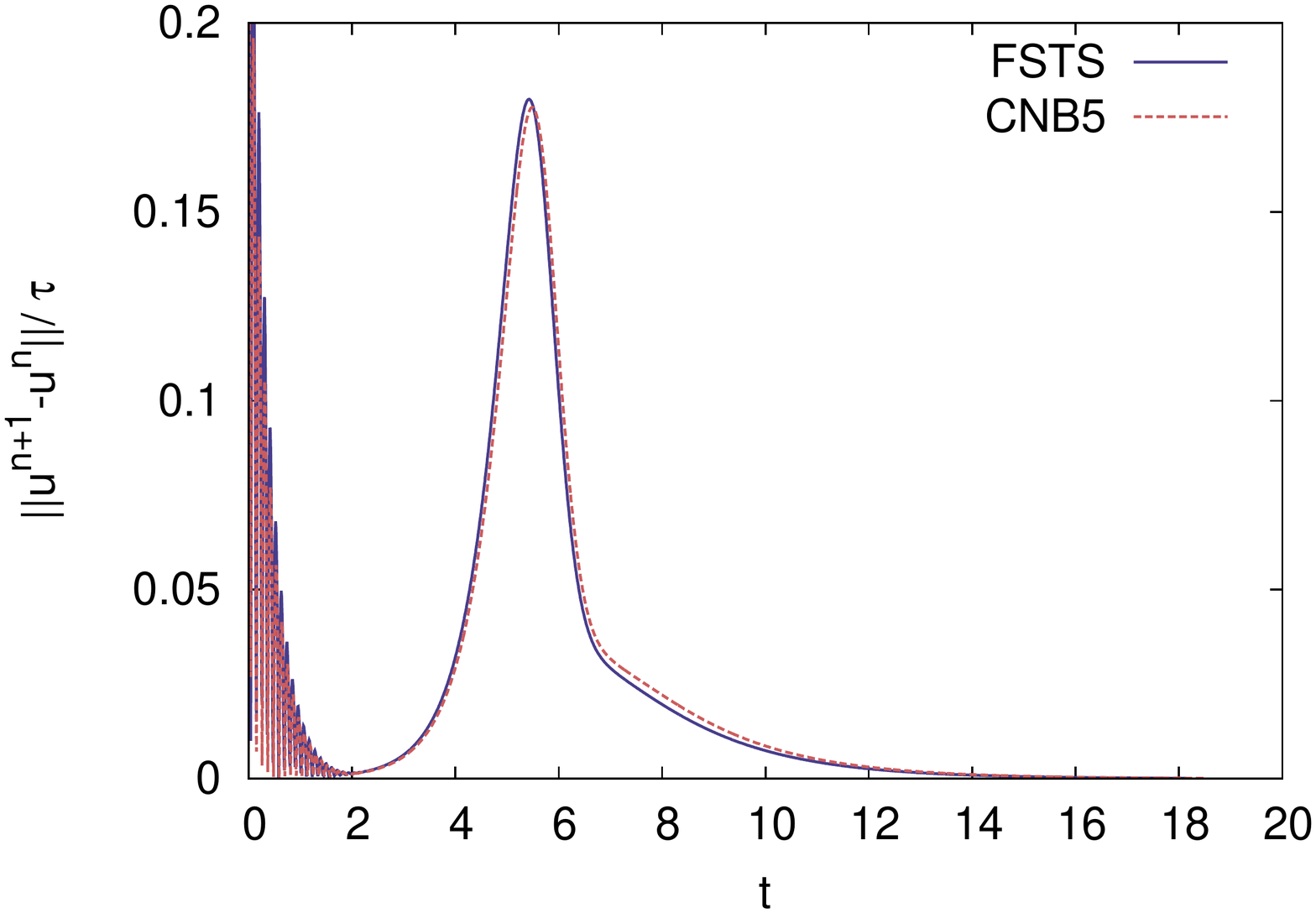}
	\caption{Comparison between the convergence history of the FSTS and the modified CNB5 method for the $u$-variable. Much better agreement is seen here compared to the unmodified CN method. Though not shown, similar improvement is seen in the $v$-variabe as well. An adaptive Newton method was employed for the fully implicit time-stepping scheme.}
	\label{tcbb}
	\end{figure}

\par The convergence history of each time-stepping method using the Newton method is shown in Fig. \ref{bct}. The initial growth of the selected mode is evident, as is the decay of the other modes,
afterwhich there is the evolution to the steady state. This behaviour is reflected in the number of non-linear iterations per timestep in Fig. \ref{iter}: as the solution is changing to a greater degree between successive timesteps during the exponential growth, more non-linear iterations are needed for convergence. There is good agreement between the solutions using the BE method
and FSTS, however there is some discrepancy with the results obtained using the CN and, in fact, the solution using the CN fails to converge in the chosen time limit. From Fig. \ref{bct} it would seem that there are modes less given to decay in the CN than in the other methods. As seen in Section \ref{sec:EOC}, the CN method is prone to oscillations. To address these issues, we shall use the CNB5 method as above where we perform BE for the first five timesteps. The results are shown in Fig. \ref{tcbb}. As expected, there is much better agreement in this case.

\begin{figure}[h]
	\centering
	\subfloat[][]{
\includegraphics[keepaspectratio=true,width=.48\textwidth]{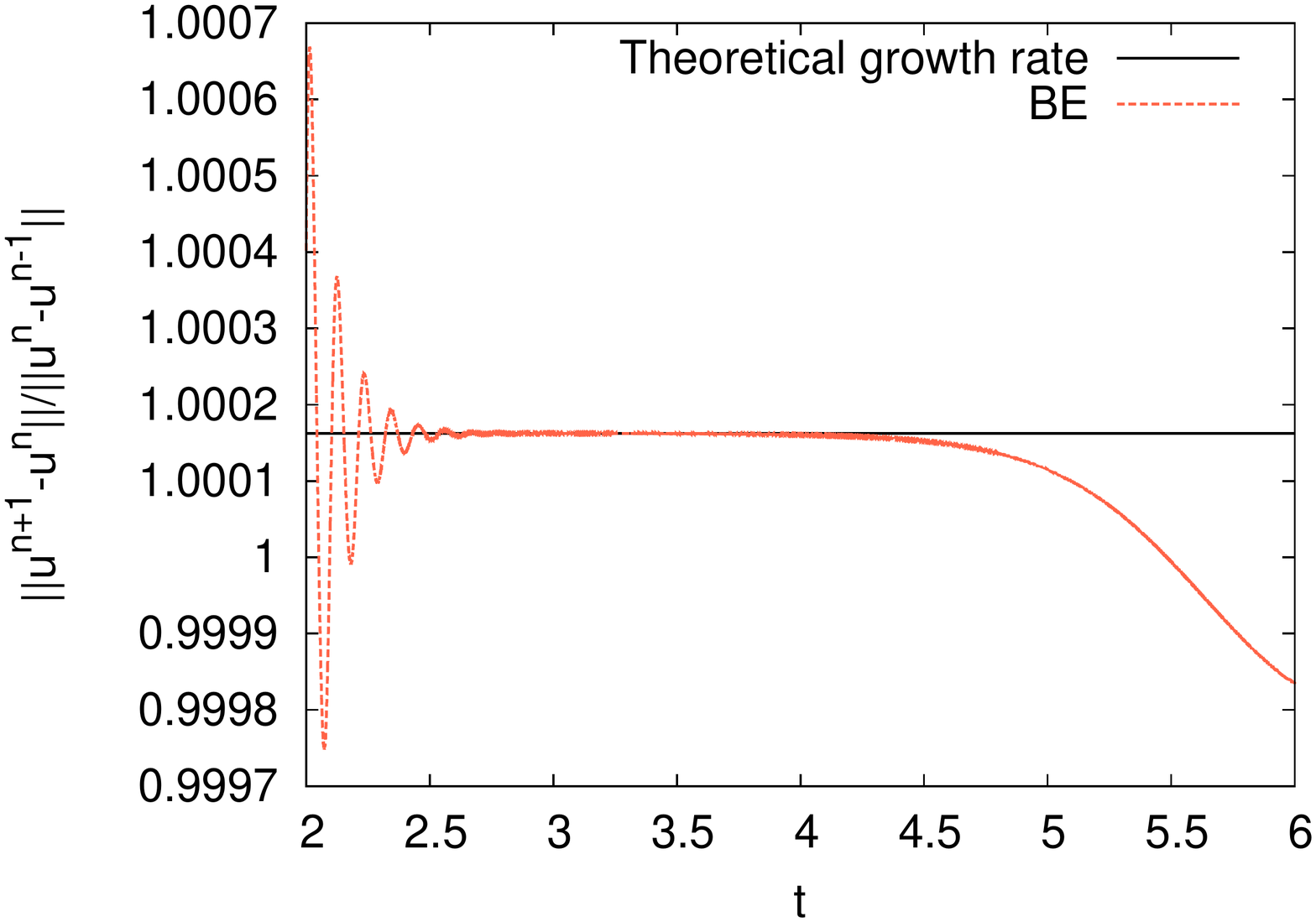}
	\label{growthratecheck:1}}
	\subfloat[][]{
\includegraphics[keepaspectratio=true,width=.48\textwidth]{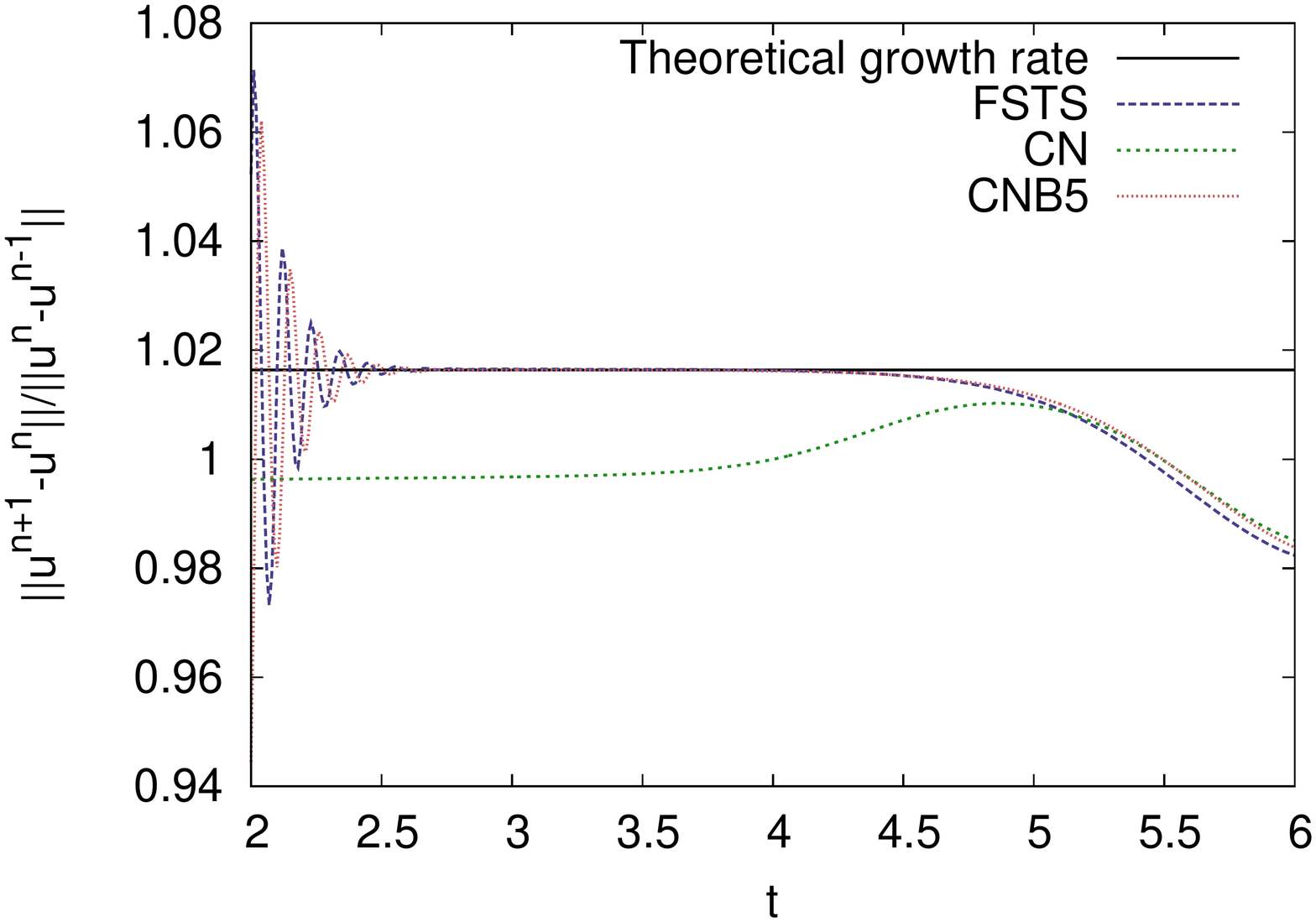}
	\label{growthratecheck:2}}
	\caption{A comparison of growth rates - the horizontal lines correspond to the theoretical growth rate as in the RHS of (\ref{theogrowthrate}). (a) BE; (b) CN, CNB5 and
FSTS. Excellent agreement during the period of mode excitation is seen between the theoretical growth rate and that obtained in
the numerical simulations using BE, CNB5 and FSTS. Using CN, this agreement is lost. An adaptive Newton method was employed for the fully implicit time-stepping scheme.}
	\label{growthratecheck}
	\end{figure}
\par The theoretical exponential growth rate of the excited mode calculated from the linear stability analysis is about $\lambda\!=\!1.6246$. So, if $u\sim e^{\lambda t}$ we should
have
\begin{equation}
\frac{||u^{n+1}-u^n||}{||u^{n}-u^{n-1}||} \sim
\frac{e^{\lambda\tau}-1}{1-e^{-\lambda\tau}}, \quad n>2.
\label{theogrowthrate}
\end{equation}
The right-hand side (RHS) is shown in the straight lines in Fig.
\ref{growthratecheck}, along with the
left-hand side (LHS) for the variable $u$ in each test. The BE and the FSTS methods show excellent agreement with the theory during the period of exponential growth, whereas the CN method shows visible deviation. The results from the modified CNB5 method, shown in Fig. \ref{growthratecheck}\subref{growthratecheck:2} show better agreement.

\section{Comparison with IMEX schemes}
\label{sec:comparison_with_imex_schemes}
\begin{figure}[thb]
	\centering	
\includegraphics[keepaspectratio=true,width=.99\textwidth]{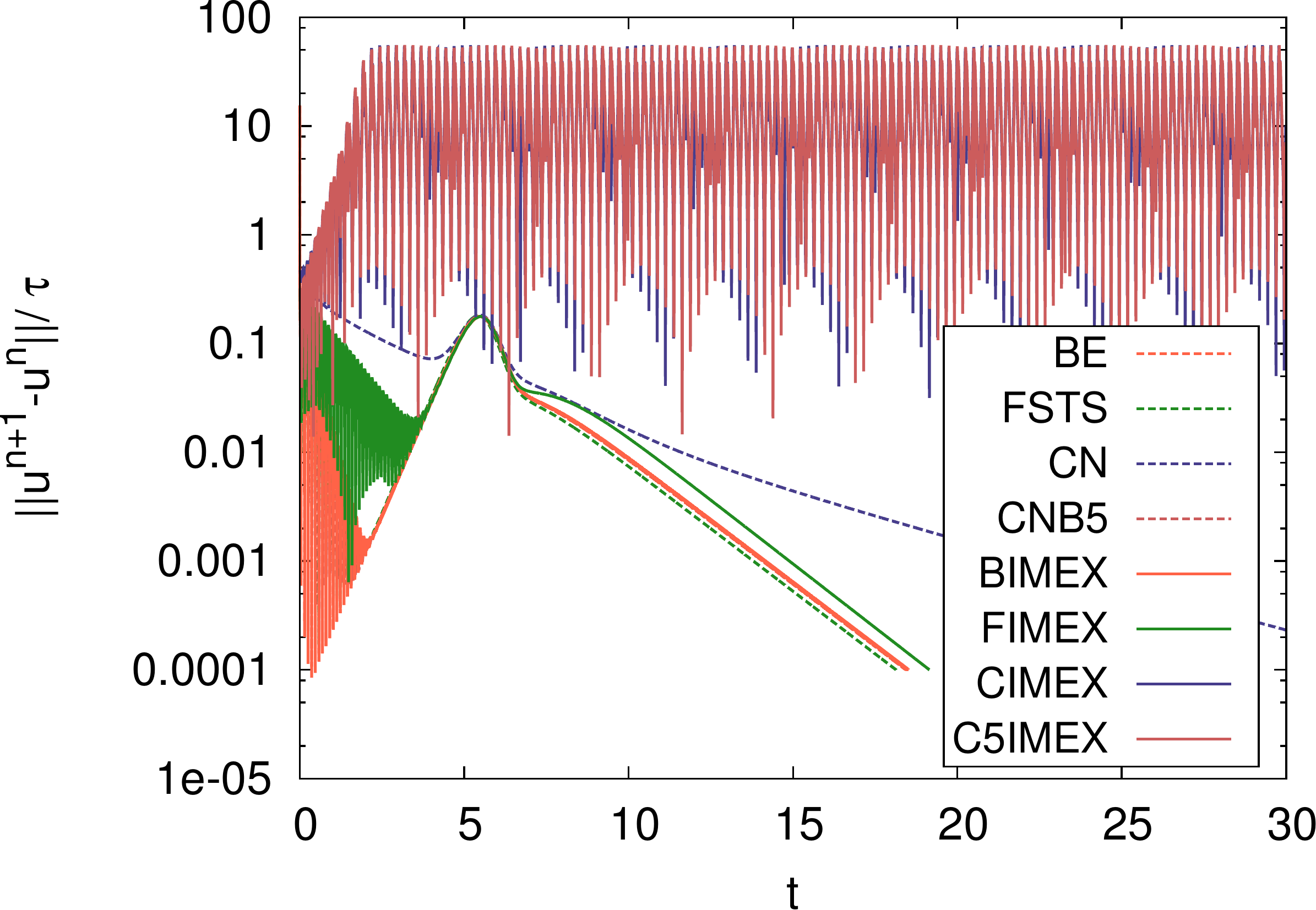}
    \caption{Comparison of the convergence history of the simulation of the Schnakenberg system for the IMEX schemes. The solid lines show the IMEX schemes and the dotted lines their fully implicit counterparts.}
	\label{IMEX}
	\end{figure}

\begin{table}[ht]
      \centering
        \begin{tabular}{ l | c c c c}
                    & BIMEX & CIMEX & CBIMEX & FIMEX \\
         \hline
         end time  & 18.5024 & 30.00 & 30.00 & 19.15 \\
         elapsed CPU time ($\times 10^4$s) & 153.7 & 16.26 & 8.82 &  7.65\\  
         
             \label{iterimex}
         \end{tabular}
         \caption{Details of the simulation of the Schnakenberg system using the IMEX schemes with one Picard iteration. Shown are the end time of the simulations and the elapsed CPU time. Note the non-convergence of CIMEX and C5IMEX.}
         \label{table:iterimex}
\end{table}
\begin{table}[ht]
      \centering
        \begin{tabular}{ l | c c c c}
                    & BE & CN & CNB5 & FSTS \\
         \hline
         end time  & 18.4571 & 30.00 & 18.49 & 18.14 \\
         elapsed CPU time ($\times 10^4$s) & 793.9 & 124.9 & 87.3  &  9.97 \\  
         
         \end{tabular}
         \caption{Details of the simulation of the Schnakenberg system using the one Newton iteration. Shown are the end time of the simulations and the elapsed CPU time.}
         \label{table:iternimex}
\end{table}

\par For the sake of comparison, the simulation was carried out without using the iterative techniques to solve the non-linearities but otherwise leaving all other parameters, outlined in Section \ref{sec:num_results}, unaltered. Instead, only one Picard iteration was performed at each timestep which is equivalent to
linearising the non-linear terms to obtain an IMEX scheme. This was done with BE (BIMEX), CN (CIMEX), CNB5 (C5IMEX) and FSTS (FIMEX). The results are  shown in Table \ref{table:iterimex} and Fig. \ref{IMEX}. The results most affected are the CN schemes, both modified and unmodified. Under the IMEX schemes, neither
converges and both suffer from colossal oscillations spanning magnitudes of order $10^2$. With FIMEX, the change is not as dramatic, however in Fig. \ref{IMEX} the convergence to the steady state towards the end is visibly different. The least affected is BN, with the BIMEX and BN lines in Fig. \ref{IMEX} indiscernible from each other. This is to be expected since the time-step $\tau$ is smaller for BN and BIMEX than for the second order methods and so the linearisation of the nonlinear terms is more accurate. The solution for CIMEX, 
C5IMEX and FIMEX are compromised due to the fact that being explicit time-stepping schemes, their region of stability is reduced \cite{Mad2006}. The value $\tau\!=\!0.01$ does not fall within this region. However, as we saw in Section \ref{sec:num_results}, this value of $\tau$ gave good results when we treated the non-linearities implicitly. This highlights the fact that in using a fully implicit scheme we can allow ourselves to use a larger timestep  thereby allowing for greater numerical stability.

To further investigate this observation, the simulations were carried out again using only one Newton iteration instead of one Picard one. The results agreed well without controlling the number of Newton iterations as in Section \ref{sec:num_results} and the convergence histories are virtually the same. The information about the end time of simulation and the CPU time taken is shown in \ref{table:iternimex}. Comparing Tables \ref{table:iter} and \ref{table:iternimex} one can see good agreement in the end time of simulation. Hence the time taken by the fully adaptive Newton and the restricted single Newton methods is similar; there does not seem to be much or any speed up  or gained by using the fully adaptive versus the single Newton iteration. This is not at all surprising since the full Newton varied between one and two iterations per timestep, so only taking one iteration does not actually significantly reduce the total elapsed CPU time required.  

\section{Applications to other geometries}\label{sec:geometries}
The methods outlined are readily applicable to complex geometries in multi-dimensions as well as on surfaces. Here we present some solutions of the Schnakenberg system in the bulk of the unit sphere, the unit cube and on the unit square for different parameter values. For fast and accurate simulations, we use only the FSTS to discretise in time, coupled with the Newton's method (with a single iteration at each timestep) to solve the nonlinearities implicitly. The simulations in this section were carried out using the software {\it deal.II} \cite{Bangerth2007}.
 \par
\begin{figure}[tb]
	\centering
  \hspace{1cm}
	\subfloat[][]{
\includegraphics[scale=.4]{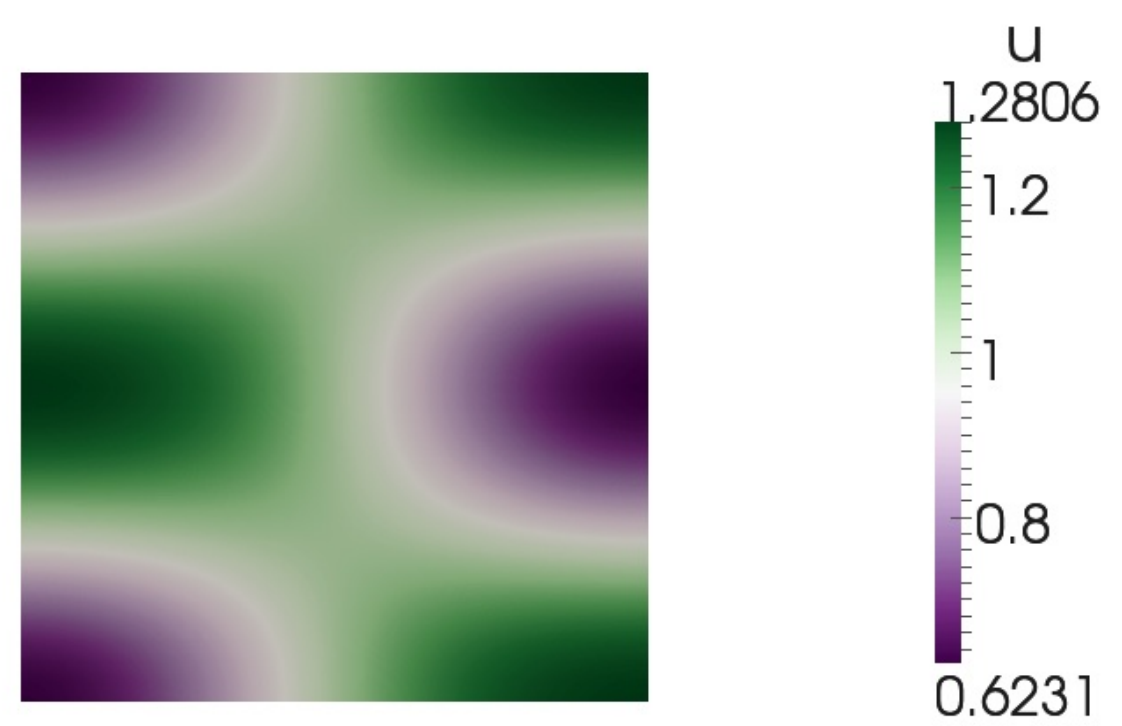}
\label{fig:schnak_square_21}}
 \hspace{3cm}
     \subfloat[][]{
\includegraphics[scale=.4]{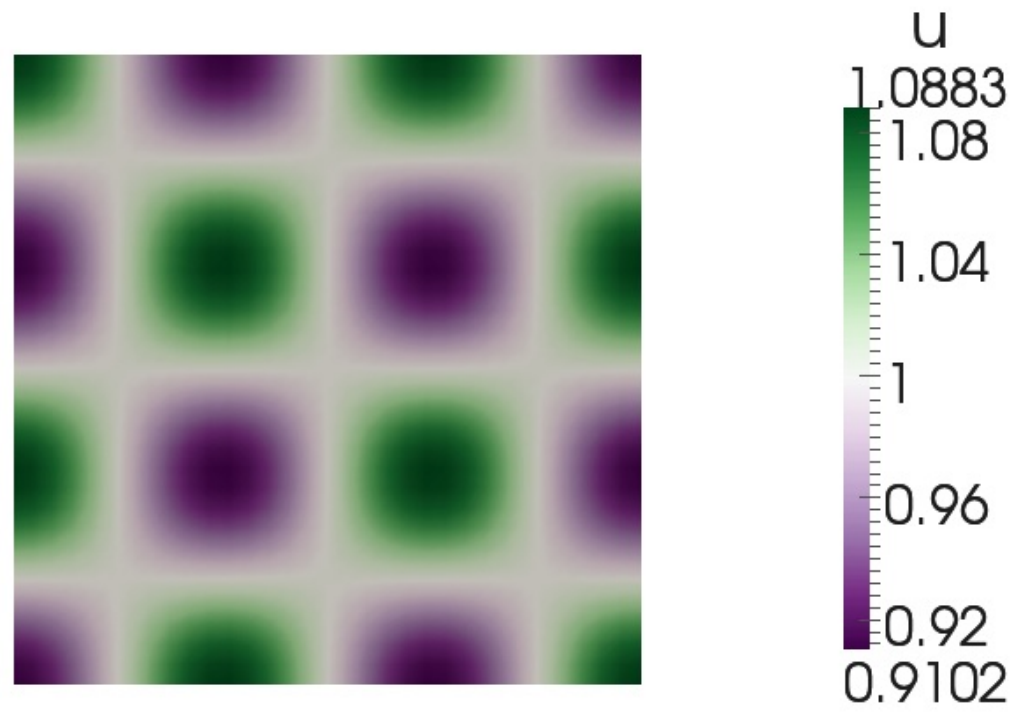}
 \label{fig:schnak_square_33}}
  \newline
    	\subfloat[][]{
\includegraphics[width=0.48\textwidth,height=3.6cm]{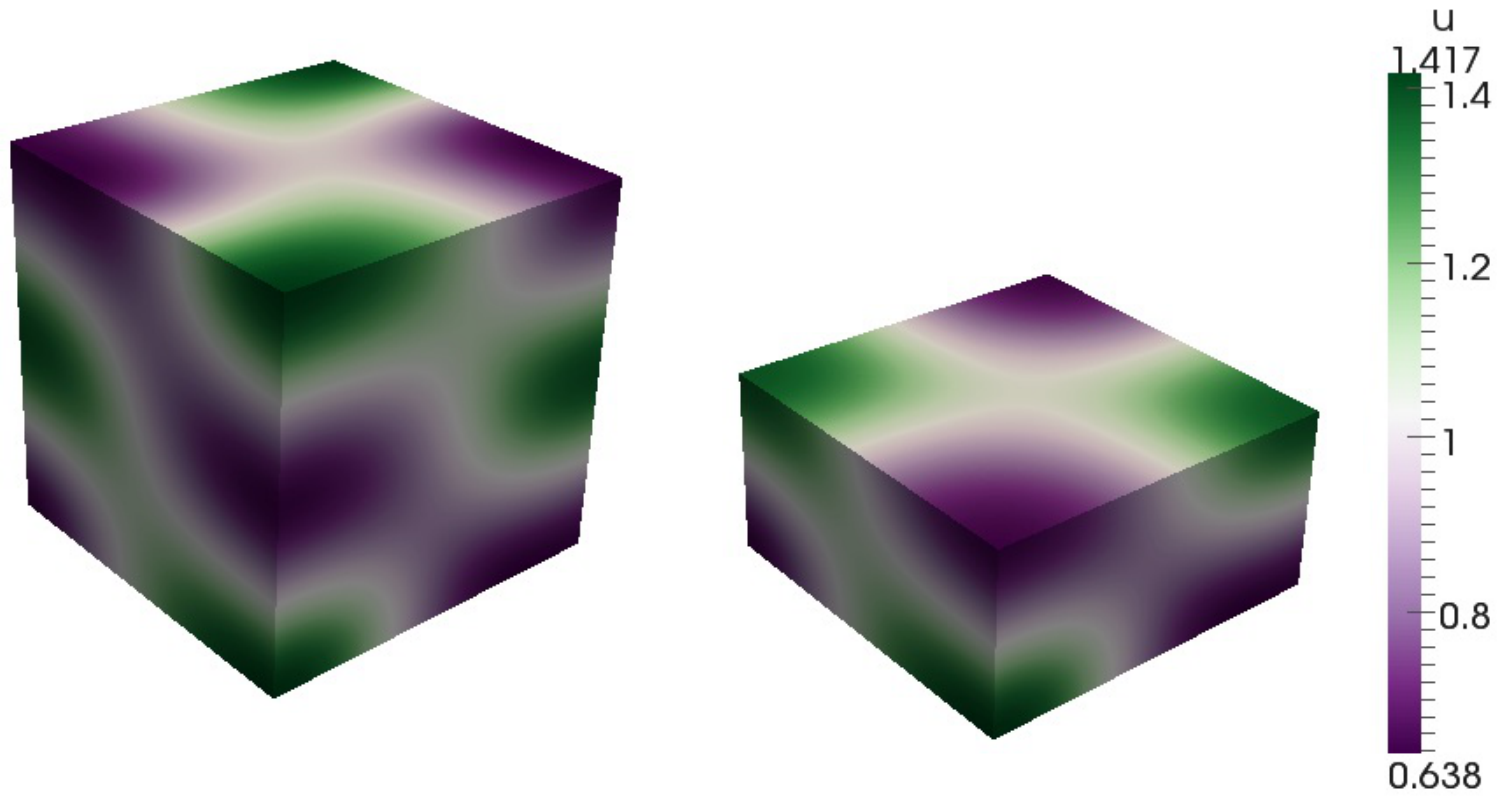}
\label{fig:schnak_cube_21}}
     \subfloat[][]{
\includegraphics[width=0.48\textwidth,height=3.6cm]{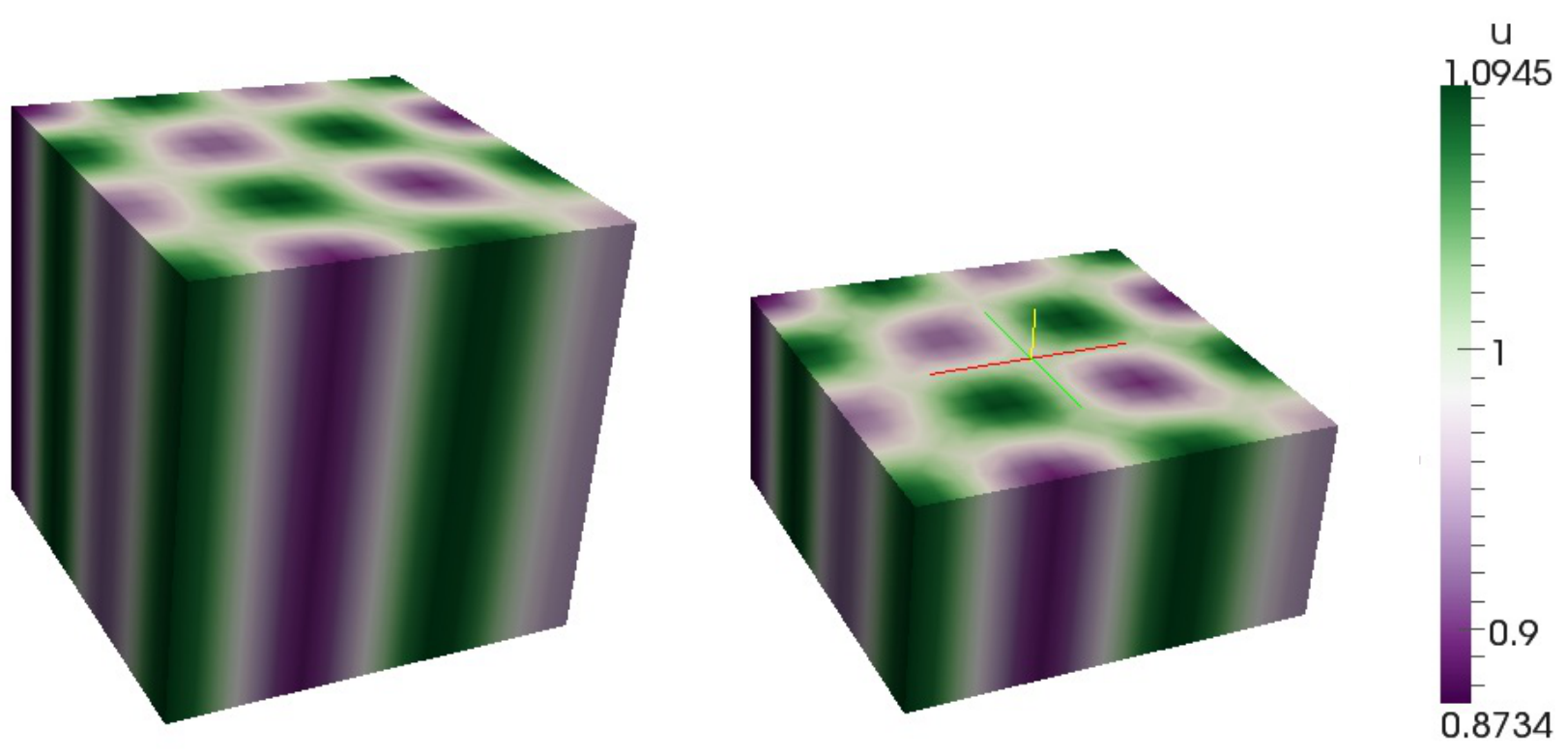}
 \label{fig:schnak_cube_33}}
 \newline
	\subfloat[][]{
\includegraphics[width=0.46\textwidth,height=3.5cm]{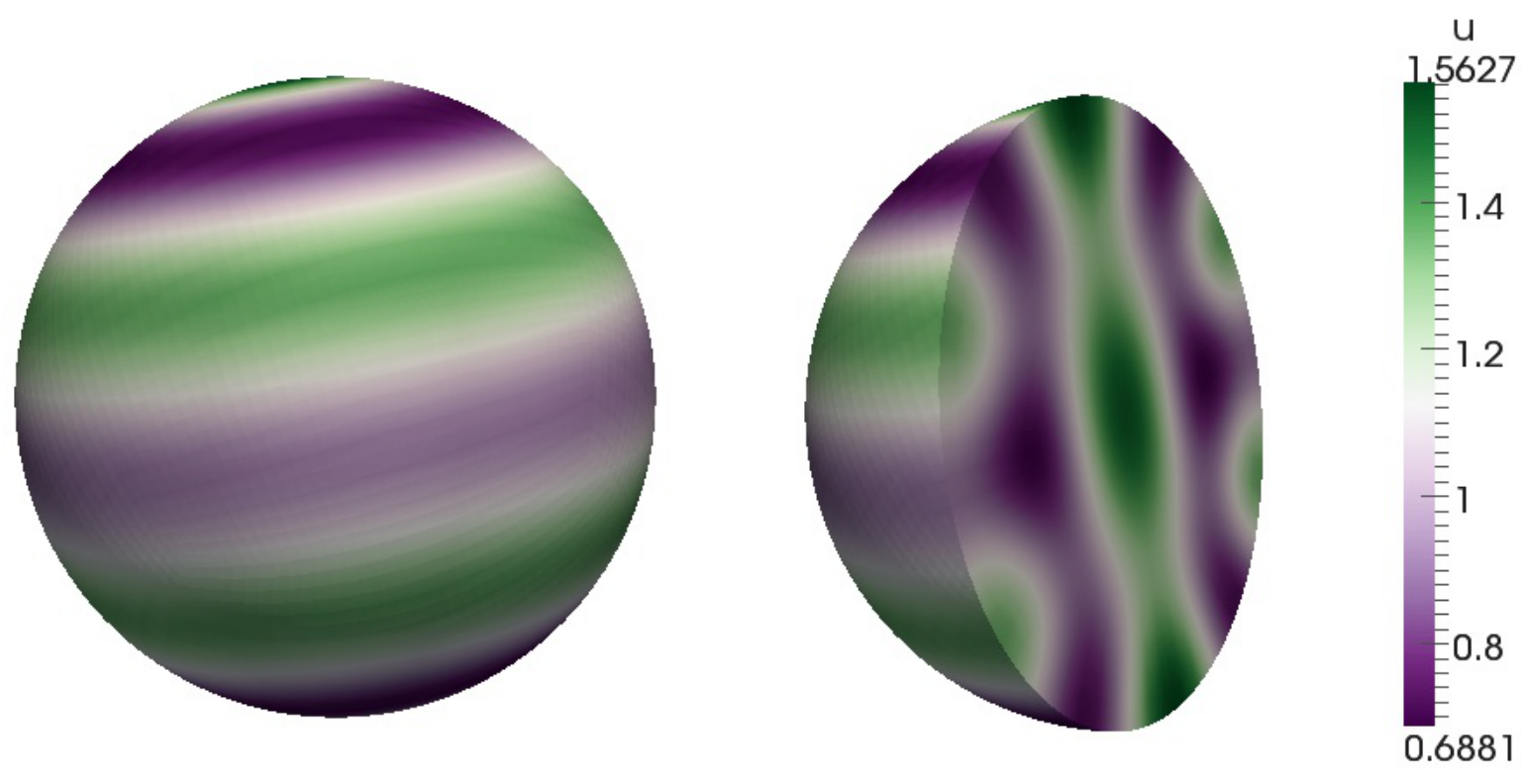}
\label{fig:schnak_sphere_21}}
     \subfloat[][]{
\includegraphics[width=0.46\textwidth,height=3.5cm]{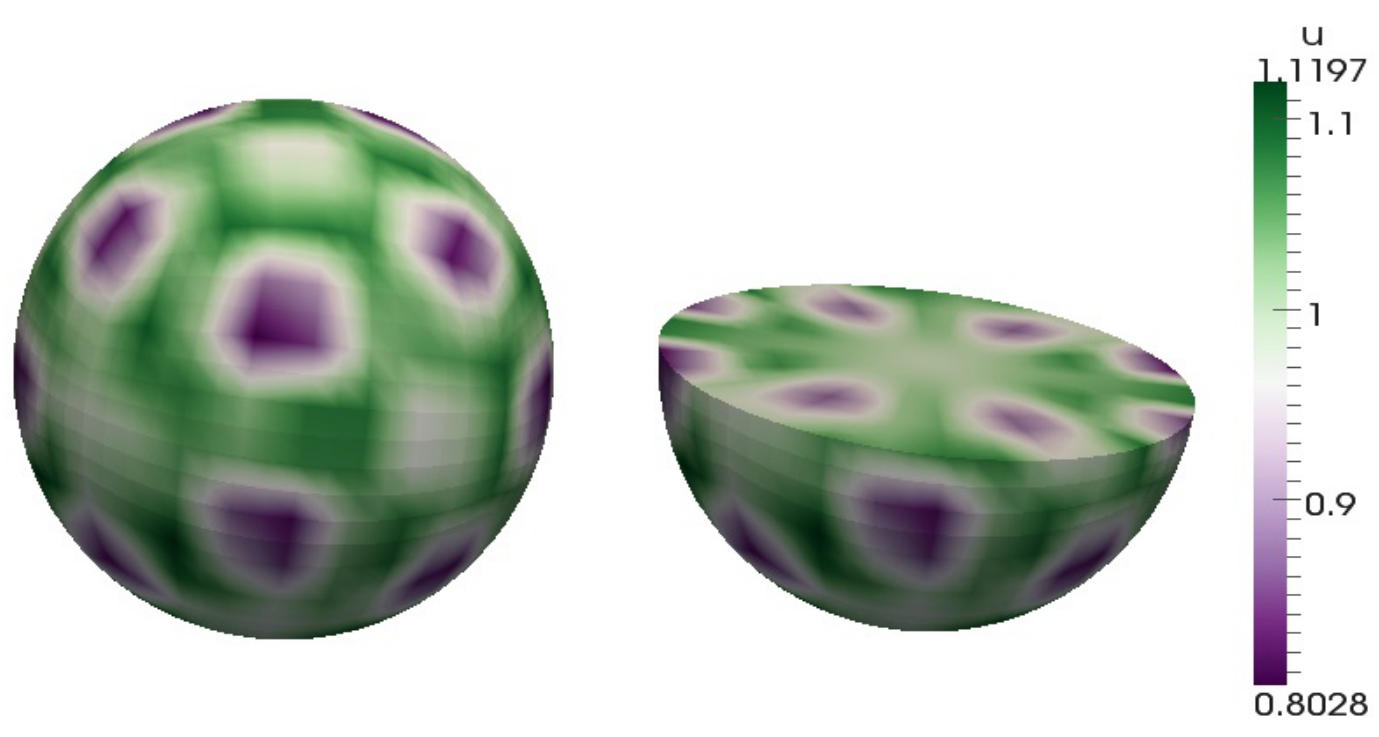}
 \label{fig:schnak_sphere_33}}
\caption{Solution for the variable $u$ of the Schnakenberg system on the square
(top row), the bulk of the cube (middle row) and the bulk of the sphere (bottom
row). Parameter values were $a\!=\!0.1$, $b\!=\!0.9$ and (a), (c), (e)
$d\!=\!9.1676$, $\gamma\!=\!176.72$ and (b), (d), (f)  $ d\!=\!8.6076$ ,
$\gamma\!=\!535.09$. Part of the domain has been cut away and shown on the right
in the 3D geometries  to reveal some internal structure. A single Newton iteration at each timestep was employed for the fully implicit time-stepping scheme.}
\label{fig:schnak_other_geoms}
\end{figure}

 We choose parameters to isolate modes on the square and use the same parameters on the 3D geometries. As detailed before, for the unit square the eigenmodes of (\ref{laplaceeigen}) have the form $\cos(n\pi x)\cos(m\pi y)$ for 
$n,m\in\mathbb{Z}$ with eigenvalues $k^2\!=\!n^2+m^2$. We choose the parameter values  $a\!=\!0.1$, $b\!=\!0.9$ and firstly $d\!=\!9.1676$, $\gamma\!=\!176.72$, then secondly $ d\!=\!8.6076$, $\gamma\!=\!535.09$. The first set isolates modes corresponding to $(n,m)\!=\!(2,1)$, whilst the second isolates those corresponding to $(3,3)$ \cite{Mad2000}. The results of such calculations are shown in Fig.  \ref{fig:schnak_other_geoms}. A wide variety of patterns are obtained in these geometries. 

\section{Conclusions}
\label{sec:conclusion}
Overall elapsed CPU times in numerical experiments shown in Table \ref{table:iter} demonstrate that the fractional $\theta$-method (FSTS) is about 130 times faster than the backward Euler and 15 times faster than the Crank-Nicholson and its modifications when both the Picard iteration and the Newton method are employed as fully implicit solvers. Furthermore, numerical tests show that a single Newton iteration at each timestep outperforms the Picard iteration. However, there is a drawback; from Table \ref{table:iter} it is clear that the Newton method is slower than the Picard iteration in treating non-linearities arising in reaction-diffusion systems when the backward Euler and the Crank-Nicholson method and its modifications are employed. The extra time needed is attributed to the need by the Newton method to use the GMRES solver to invert the non-symmetric systems, whereas the Picard iteration uses the CG solver for the symmetric systems which is quicker.  In the FSTS method the non-linear system \eqref{schnakfefstseta2} is slightly different than the others and takes less iterations for the GMRES solver to solve.
\par The BE scheme and FSTS agreed well with each other, however the CN method had to be modified in order to gain such an agreement. The modified CN method performs very well; it is second order, is more straightforward to implement than the FSTS and is unconditionally stable. However,  even though the FSTS does require two extra equations to be solved per timestep, the extra equations are linear and as such do not add any significant amount of work to be done to calculate the solution. As a result, we strongly recommend using the FSTS coupled with the Newton method (with a single iteration at each timestep) as the preferred fully implicit time-stepping scheme for reaction-diffusion systems on stationary and sometimes continuously evolving domains, volumes and surfaces. 
\par The comparison between the IMEX schemes and the fully implicit schemes demonstrate how linearisation affects the results. Solving non-linearities implicitly alleviated the need for smaller timesteps. However, this alleviation came at the price of solving equations that were non-linear. The choice between IMEX and fully implicit schemes then lies in the trade-off between the size of
the timestep and the number of extra non-linear iterations at each timestep. 
\par It is to be expected that such results will also apply to other reaction-diffusion systems such as the Gierer-Meinhardt model \cite{GM1972}, the Thomas model \cite{TK1976} and Brusselator model \cite{PL1968}. In the case of different boundary conditions, the inhomogeneous steady state and the Turing space may change but the same method of discretisation may be used, so that the results presented in this article are also applicable.

Although we have presented results for stationary domains, it is natural to extend the numerical analysis of implicit solvers to domains and surfaces that change continuously, i.e. growing domains and evolving surfaces. It is when solving non-autonomous systems of PDEs posed on evolving domains and surfaces (as well as in the bulk) that a single Newton iteration will become crucially
important both as an accurate solver as well as an aid for large  computational savings. The theory of PDEs posed on evolving domains and surfaces (as well as in the bulk) is a fast burgeoning research area with many applications in cell motility, plant biology and developmental biology where such models are
routinely used  \cite{Elliott2012,Neilson2011,Venkataraman2011,Lakkis2013}.

\Appendix
\section{Convergence of the Picard iteration}\newline
Suppose we have a system of PDE's of the form
\begin{equation} 
 \left\{
 \begin{aligned}
 & \frac{ \partial u}{\partial t} - d_1\nabla^2u  =  c_1+\gamma f(u,v)u, \\
 & \frac{ \partial v}{\partial t} - d_2\nabla^2v =  c_2 + \gamma g(u,v)v, \indom{\Omega}\,,\mbox{for } t>0\,,
 \end{aligned}
 \right.  \label{geneqn}
 \end{equation}
for the diffusion coefficients $d_1$, $d_2\!>\!0$, constant source terms $c_1$, $c_2\!\in\!\mathbb{R}$ and functions $f\!:\!\Omega\!\rightarrow\!\mathbb{R}$ and $g\!:\!\Omega\!\rightarrow\!\mathbb{R}$ which are in general non-linear, but globally Lipshitz continuous with Lipshitz constants $L_f$ and $L_g$ respectively, and with no explicit dependence on $x$ or $t$. Even though the Schnakenberg non-linearities may not be globally Lipshitz in general, the solution remains in an invariant set and therefore a Lipschitz constant may still be found in this region \cite{Venkataraman2011}. We further suppose that $u$ and $v$ always remain positive and that $f$ and $g$ are bounded above by some constants $\beta_f$ and $\beta_g$ respectively. Let $T_h$ be a quasi-uniform triangulation of $\Omega$, and $S_h$ be the space of continuous functions on $T_h$ that are piecewise linear on each element. Recasting into the weak form and using the BE method to discretise time as above we solve for $(u^{n+1},v^{n+1})\!\in\!S_h\!\times\!S_h$ from the system
 \begin{equation}  
   \left\{
 \begin{aligned}
 & (\,\frac{u^{n+1}-u^n}{\tau},\phi_h)+d_1(\nabla u^{n+1},\nabla\phi_h)= (c_1,\phi_h)+ (f(u^{n+1},v^{n+1})u^{n+1},\phi_h), \\
 & (\,\frac{v^{n+1}-v^n}{\tau},\phi_h)+d_2(\nabla v^{n+1},\nabla\phi_h)= (c_2,\phi_h)+ (g(u^{n+1},v^{n+1})v^{n+1},\phi_h), \quad \forall \phi_h\in S_h,
 \end{aligned}
 \right.  \label{geneqn:fem}
 \end{equation}
where $(\cdot\,,\cdot)$ denotes the inner product in $L_2$. We would like to use Picard iteration to solve the non-linear equation iteratively, each iteration requiring a linear problem to be solved. The $(k+1)$-th iterate is then
 \begin{equation}  
   \left\{
 \begin{aligned}
 & \left(\frac{u^{n+1}_{k+1}-u^n}{\tau},\phi_h \right)+d_1(\nabla u^{n+1}_{k+1},\nabla\phi_h )= (c_1,\phi_h) + \left(f(u^{n+1}_k,v^{n+1}_k)u^{n+1}_{k+1},\phi_h\right), \\
 & \left(\frac{v^{n+1}_{k+1}-v^n}{\tau},\phi_h \right)+d_2(\nabla v^{n+1}_{k+1},\nabla\phi_h)= (c_2,\phi_h) + \left(g(u^{n+1}_k,v^{n+1}_k)v^{n+1}_{k+1},\phi_h\right).
 \end{aligned}
 \right.  \label{geneqn:pic}
 \end{equation}
Subtracting (\ref{geneqn:fem}) from (\ref{geneqn:pic}) and choosing appropriate $\phi_h$ we then have
 \begin{equation}  
   \left\{
 \begin{aligned}
 & \frac{1}{\tau}||u_{k+1}-u^{n+1}||^2 \leq \left( f(u_k,v_k)u_{k+1}-f(u^{n+1},v^{n+1})u^{n+1},u_{k+1}-u^{n+1}\right), \\
 & \frac{1}{\tau}||v_{k+1}-v^{n+1}||^2 \leq \tau\left( g(u_k,v_k)v_{k+1}-g(u^{n+1},v^{n+1})v^{n+1},u_{k+1}-u^{n+1}\right),
 \end{aligned}
 \right.  \label{geneqn:pic1}
 \end{equation}
 where $||\cdot||$ denotes the $L_2$-norm. We have dropped the superscript $(n+1)$ from the Picard iterates in (\ref{geneqn:pic1}) for clarity. Concentrating on the equation in the $u$-variable, we rewrite the  RHS as
 \begin{equation}
  \begin{split}
 &\Big(\big(f(u_k,v_k)-f(u^{n+1},v^{n+1})\big)(u_{k+1}-u^{n+1}),u_{k+1}-u^{n+1}\Big) \\
 & \quad +2\Big(\big(f(u_k,v_k)-f(u^{n+1},v^{n+1})\big)u^{n+1},u_{k+1}-u^{n+1}\Big) \\
 & \quad +\left( f(u^{n+1},v^{n+1})u_{k+1}-f(u_k,v_k)u^{n+1},u_{k+1}-u^{n+1}\right),
 \end{split}
 \end{equation}
 and obtain
 \begin{equation}
 \begin{split}
 & \frac{1}{\tau} ||u_{k+1}-u^{n+1}||^2\\
 & \quad  \leq ||f(u_k,v_k)-f(u^{n+1},v^{n+1})||\Big(
          \, ||(u_{k+1}-u^{n+1})^2||
        +2||u^{n+1}||_{\infty}||u_{k+1}-u^{n+1}||\Big) \\
  &  \qquad +\beta_f ||u_{k+1}-u^{n+1}||^2 \\
  &  \quad \leq L_f\big(C_1+2||u^{n+1}||_{\infty}\big) ||u_{k+1}-u^{n+1}||\,||\bs{\xi}_k-\bs{\xi}^{n+1}||\ + \beta_f ||u_{k+1}-u^{n+1}||^2,
 \end{split}
 \end{equation}
 for large enough $C_1$, and where $\bs{\xi}_k\!=\!(u_k,v_k)$, and similarly for $\bs{\xi}^{n+1}$. Upon rearrangement, this becomes
 \begin{equation}  
  ||u_{k+1}-u^{n+1}|| \leq L_f\frac{\tau}{1-\tau\beta_f}\big(C_1+2||u^{n+1}||_{\infty}\big) ||\bs{\xi}_k-\bs{\xi}^{n+1}||.
   \label{geneqn:u_estim}
 \end{equation}
Similarly,
 \begin{equation}  
  ||v_{k+1}-v^{n+1}|| \leq L_g\frac{\tau}{1-\tau\beta_g}\big(C_2+2||v^{n+1}||_{\infty}\big) ||\bs{\xi}_k-\bs{\xi}^{n+1}||.
   \label{geneqn:v_estim}
 \end{equation}
Since $||\bs{\xi}_{k+1}-\bs{\xi}^{n+1}||\!\leq\!||u_{k+1}-u^{n+1}||\!+\!||v_{k+1}-v^{n+1}||$
we obtain the estimate
 \begin{equation}  
  ||\bs{\xi}_{k+1}-\bs{\xi}^{n+1}|| \leq L\frac{\tau}{1-\tau\beta}\big(
                   C+2||\bs{\xi}^{n+1}||_{\infty}\big)
                    ||\bs{\xi}_k-\bs{\xi}^{n+1}||,
   \label{geneqn:xi_estim}
 \end{equation}
 where $L\!=\!\max\{L_f,L_g\}$, $\beta\!=\!\max\{\beta_f,\beta_g\}$, $C\!=\!\max\{C_1,C_2\}$. Let $\bs{\Xi}$ be the exact solution $(U,V)$ to the original problem (\ref{geneqn}). Then
 \begin{equation}
 ||\bs{\xi}^{n+1}||_{\infty}\leq ||\bs{\Xi}(t^{n+1})-\bs{\xi}^{n+1}||_{\infty}+ ||\bs{\Xi}(t^{n+1})||_{\infty} \leq C^{\prime}(\bs{\Xi})(h^2+\tau) + ||\bs{\Xi}(t^{n+1})||_{\infty},
 \end{equation}
for some constant $C^{\prime}(\bs{\Xi})$. Thus, redefining constants as appropriate, we finally obtain
 \begin{equation}  
 \begin{split}
  ||\bs{\xi}_{k+1}-\bs{\xi}^{n+1}|| & \leq CL\frac{\tau}{1-\tau\beta}\big(
                   1+C^{\prime}(\tau+h^2)\big)
                    ||\bs{\xi}_k-\bs{\xi}^{n+1}|| \\
  &\quad \colonequals \lambda\, ||\bs{\xi}_k-\bs{\xi}^{n+1}||.
   \label{geneqn:final_estim}
   \end{split}
\end{equation}
The coefficient $\lambda$ will be less than $1$ for small enough $\tau$ and $h$, thus showing that the method converges. In particular, $\lambda\rightarrow 0$ as $\tau\rightarrow 0$ and $h\rightarrow 0$ simultaneously.
\newline

\noindent{ \textnormal{\textbf{Acknowledgements.}
This work (AM) is partly supported by the following grants: the Engineering and Physical Sciences Research Council (EP/J016780/1), the London Mathematical Society (R4P2) and  the British Council through its {\it UK-US New Partnership Fund (PMI2)}.  AHC was supported partly by the University of Sussex and partly by the Medical Research Council DTA case award.}}
\par \mbox{}
\bibliography{references}
\bibliographystyle{plain}

\end{document}